\documentclass{article}
\pdfoutput=1

\usepackage[pdftex]{graphicx}
\usepackage{xy}
\xyoption{all}
\usepackage{amssymb, amsmath, amscd}
\usepackage{amsthm}
\usepackage{mathrsfs}

\usepackage{blkarray}

\usepackage{pgfplots}
\usepackage{tikz}
\usetikzlibrary{positioning}
\usetikzlibrary{calc}

\usetikzlibrary{arrows,arrows.meta}

\newcommand*{\Cite}{\cite}

\newcommand{\Ref}[1]{\textit{\ref{#1}}}
\newcommand{\RefPar}[1]{\textbf{(\ref{#1})}}

\def\N{\mathbf{N}}
\def\R{\mathbf{R}}
\def\Q{\mathbf{Q}}
\def\Z{\mathbf{Z}}
\def\C{\mathbf{C}}

\def\dAnn{\mathbf{Rng}^{\partial}}

\newcommand{\A}[2]{\mathbf{A}_{#2}^{#1}}
\newcommand{\Pn}[2]{\mathbf{P}_{#2}^{#1}}

\newlength\longueurdelaboite
\newcommand{\boite}[1]
{
\settowidth{\longueurdelaboite}{#1}
\begin{minipage}{\longueurdelaboite}
#1
\end{minipage}
}

\newlength\longueurdudiagramme 
\newcommand{\fonctionb}[4]
{
\settowidth{\longueurdudiagramme}{$\xymatrix@R=0mm{
#1 \ar[r] & #2 \\
#3 \ar@{|->}[r] & #4
}$} 
\begin{minipage}{\longueurdudiagramme}
$\xymatrix@R=0mm{
#1 \ar[r] & #2 \\
#3 \ar@{|->}[r] & #4
}$
\end{minipage}
}

\newcounter{compteur}
\renewcommand{\thecompteur}{\alph{compteur}}
\newcommand{\point}{\smallskip \refstepcounter{compteur} \textit{(\thecompteur)}{} }

\newcommand{\remark}[1]{\noindent  \textit{Remark.\,---\!-- } #1 $\diamondsuit$}

\newcommand{\remarks}[1]{\noindent\setcounter{compteur}{0}\textit{Remarks.\,---\!-- }\point #1 $\diamondsuit$}

\newcommand{\example}[1]{\noindent  \textit{Example.\,---\!-- } #1 $\triangle$}

\newcommand{\examples}[1]{\noindent   \setcounter{compteur}{0}   \textit{Examples.\,---\!-- }\point  #1 $\triangle$}

\newcommand{\pproof}[1]{\noindent  \textit{Proof.\,---\!-- } #1$\blacksquare$}

\newcommand{\pproofbis}[2]{\noindent  \textit{Proof {#1}.\,---\!--} #2 $\blacksquare$}

\DeclareMathOperator{\colim}{{\mathrm{colim}}}

\DeclareMathOperator{\limind}
{%
\underrightarrow{
\raisebox{-1pt}[0pt][1.5pt]{}
\colim
\raisebox{-1pt}[0pt][1.5pt]{}
}}

\def\impl{\Longrightarrow}
\def\longto{\mathop{\longrightarrow}}
\def\tq{\mid}
\def\Sp{\mathrm{Spec\,}}
\def\fa{\forall}
\def\xt{\exists}
\newcommand{\ie}{\emph{ie} }
\newcommand{\aand}{\text{ and }}

\def\ssi{\iff}
\def\Id{\mathrm{Id}}

\newcommand{\ob}[1]{\mathrm{ob}\left ( {#1} \right )}

\def\dSp{\mathrm{Spec}^{\partial}}
\def\diffSp{\mathrm{diff\text{-}Spec\,}}
\newcommand{\dSch}{\mathbf{Sch}^{\partial}}
\newcommand{\Traj}[1]{\text{Traj}_{\vec{\mathscr{#1}}}}
\newcommand{\PrSh}[1]{\mathbf{PrSh}(#1)}
\newcommand{\Sh}[1]{\mathbf{Sh}(#1)}
\newcommand{\AbPrSh}[1]{\mathbf{PrSh}_{\textit{Ab}}(#1)}
\newcommand{\RngPrSh}[1]{\mathbf{PrSh}_{\textit{Rng}}(#1)}
\newcommand{\dRngPrSh}[1]{\mathbf{PrSh}_{\textit{Rng}^\partial}(#1)}
\newcommand{\AbSh}[1]{\mathbf{Sh}_{\textit{Ab}}(#1)}
\newcommand{\RngSh}[1]{\mathbf{Sh}_{\textit{Rng}}(#1)}
\newcommand{\dRngSh}[1]{\mathbf{Sh}_{\textit{Rng}^\partial}(#1)}

\newtheorem{theorem}{Theorem}[section] 
\newtheorem{lemma}[theorem]{Lemma}
\newtheorem{corollary}[theorem]{Corollary}
\newtheorem{definition}[theorem]{Definition}
\newtheorem{proposition}[theorem]{Proposition}

\newtheorem*{propositionintro}{Proposition} 
\newtheorem*{theoremeintro}{Theorem}
\newtheorem{notation}[theorem]{Notation}
\newtheorem*{lemmaintro}{Lemma}

\makeatletter
\def\footnoterule{\kern-.4\p@
        \hrule\@width 5pc\kern10\p@\kern-\footnotesep}
\def\@makefnmark{\hbox{$\,\m@th^{(\@thefnmark)}$}}
\def\@makefntext{\parindent0pt\sloppy\indent\@makefnmark}
\makeatother

\title{
{Vector fields and differential schemes}}

\author{%
Colas Bardavid \\
IRMAR --- UMR 6625 du CNRS\\
Universit\'e de Rennes 1 Campus de Beaulieu \\
35042 Rennes CEDEX FRANCE
}

\begin{document}

\maketitle
 \begin{center}  
\rule{10 cm}{.5pt} \end{center}
\begin{small}
\begin{center}
\begin{minipage}{10cm}
\textbf{Abstract --} For a scheme $X$ and $\vec{\mathscr{V}}$ a vector field on $X$, we define the leaves of $\vec{\mathscr{V}}$ (in a purely algebraic context). Given $x\in X$, we prove that there exists a smallest leaf $\eta$ containing $x$:  we call it the trajectory of $x$ and establish some useful properties for it. We take this point of view to give a geometrical interpretation of other works about differential schemes. Our main result is this natural property: it is always possible to extend, in a unique way,  a constant section defined over $U$ to the open set $U^\delta$ generated by $U$ under the action of $\vec{\mathscr{V}}$. We use these techniques to compare three sheaves that have been defined over the differential spectrum.
\end{minipage}
\end{center}
\end{small}
 \begin{center}  
\rule{10 cm}{.5pt} \end{center}

\begin{center}
 \begin{minipage}{10cm}
\end{minipage}
\end{center}

\newpage
\tableofcontents

\newpage

\section{Introduction}

\subsection{Differential schemes}

Differential schemes have been introduced
for the first time by William Keigher in \Cite{KeigherPremArticle}. The goal was to give to algebro-differential geometry\footnote{See the article of Buium and Cassidy in \Cite{OeuvresCompletesKolchin} for an excellent survey on algebro-differential geometry. Algebro-differential varieties appear naturally in, for instance, para\-metrized Galois theory:  in this setting,  the Galois group is not an algebraic group but a differential-algebraic group. See for instance \Cite{SingerCassidy}, \Cite{DiffAlgGroupsCass}, \Cite{Landesman},  \Cite{SingerHardouin}.} solid foundations, just as scheme theory for algebraic geometry. Despite contributions to this task by A. Buium \Cite{BuiumDiffSpec}, G. Carr\`a Ferro \Cite{CarraFerro2, KolchinSchemes} and J. Kovavic \Cite{KovacicDiffSchemes,KovacicGaloisDiffSchemes1,
KovacicGaloisDiffSchemes2,KovacicGlobalSections}, the category of differential schemes that was looked for is still missing. We  explain now quickly the several attempts made and why they failed. 
\medskip

We start with a differential ring $(A, \partial)$ and define the \emph{differential spectrum of $A$} to be 
\[
   \diffSp A := \Bigl \{\mathfrak{p} \tq \mathfrak{p}\text{ is a prime and differential ideal of $A$}\Bigr \}.
\]
The Zariski topology on $\Sp A$ induces a  topology on $\diffSp A$ called  the \emph{Kolchin topology}. The structure sheaf $\mathcal{O}_{\Sp A}$ on $\Sp A$ can also be restricted to $\diffSp A$. The restricted ringed space is called an \emph{affine differential schemes}. They are  differentially ringed spaces, whose stalks are local rings with a maximal ideal which is also differential. The category of differential schemes is defined to be the category of such differentially ringed spaces locally isomorphic to affine objects. This approach follows exaclty the approach of schemes.

\medskip

The first issue with this category is that we don't know if it has fibered products. Kovacic has studied this question (see section 14 of \Cite{KovacicDiffSchemes}) and proved that they exist for a restricted class of differential schemes (known as AAD differential schemes). The general answer to this question is unknown.
The second issue is about global sections. Unlike schemes, the natural morphism
\[A \longto \hat{A}:= \Gamma(\diffSp A, \mathcal{O}_{\diffSp A})  \]
is in general neither injective nor surjective. Under some assumptions (see Theorem 2.6 of \Cite{BuiumDiffSpec},  Theorem 10.6 of \Cite{KovacicDiffSchemes} or Theorem 8 of \Cite{TrushinPreprint}),  $\hat{A} \longto \hat{\hat{A}}$ is an isomorphism. 

\medskip

There exists another definition for differential schemes, given by Carr\`a Ferro in \Cite{KolchinSchemes} and based on another structure sheaf $\mathcal{O}_{\diffSp A}^{\text{(CF)}}$ on $\diffSp A$. It  looks more complicated than the previous one but  verifies
\[ A \cong \Gamma\Big (\diffSp A, \mathcal{O}_{\diffSp A}^{\text{(CF)}} \Big ). \]
The definition of this sheaf is based on the following lemma:
\begin{lemmaintro}[Lemma 1.5 of \Cite{KolchinSchemes}]
Let $A$ be a differential scheme. For each open subset $U$ of $\diffSp A$, there exists an open subset $U_\Delta$ of $\Sp A$ such that:
\begin{itemize}\setlength{\itemsep}{0mm}
\item[(i)] $U_\Delta \cap \diffSp A = U$;
\item[(ii)] $U_\Delta \supset \Sp A \setminus \overline{\diffSp A}$;
\item[(iii)] If $V$ is an open set of $\Sp A$ such that $V \cap \diffSp A=U$, then $V\subset U_\Delta$.
\end{itemize}
\end{lemmaintro}
Then, $\mathcal{O}_{\diffSp A}^{\text{(CF)}}(U)$ is defined to be $\mathcal{O}_{\Sp A}(U_\Delta)$. Our paper will cast more clarity to this construction and explain what are the relations between these two sheaves.

   \bigskip
   \subsection{Content of the paper}
The goal of this paper is to bring a new point of view on differential schemes and to apply it to already existing constructions. The first idea is to consider schemes instead of differential rings. In this setting, derivations are replaced by vector fields. Schemes with vector field have already been introduced and studied, in particular by Buium in \Cite{buium} and \Cite{BuiumLNM2} and Dyckerhoff in \Cite{dyckerhoff} but also by Umemura in \Cite{Umemura1}. Of course, vector fields are found in \Cite{EGA44}, but Grothendieck does not study them extensively. He leaves this study to other mathematicians\footnote{At the very beginning of \Cite{EGA44}, Grothendieck writes \emph{``Dans ce paragraphe, nous pr\'esentons, sous forme globale, quelques notions de
calcul diff\'erentiel particuli\`erement utiles en G\'eom\'etrie alg\'ebrique. Nous passons sous
silence de nombreux d\'eveloppements, classiques en G\'eom\'etrie diff\'erentielle (connexions,
transformations infinit\'esimales associ\'ees \`a un champ de vecteurs, jets, etc.), bien que ces
notions s'\'ecrivent de fa\c{c}on particuli\`erement naturelle dans le cadre des sch\'emas.''}} and that's exactly what we do in this paper.

Given a scheme $X$ with a vector field $\vec{\mathscr{V}}$, we introduce the \emph{leaves of $X$ for  $\vec{\mathscr{V}}$}. It is the points $\eta$ of $X$ that are invariant under the vector field --- \emph{ie}, the irreducible closed subsets tangent to $\vec{\mathscr{V}}$. Then, we prove that given $x \in X$, there exists a smallest leaf going through $x$, called  the \emph{trajectory of $x$ under  $\vec{\mathscr{V}}$}, denoted  by $\Traj{V}(x)$. The map $\Traj{V} : X \longto X$ satisfies natural properties allowing us to define a new topology on $X$. The open sets for this topology are the Zariski open sets of $X$ that are invariant under $\vec{\mathscr{V}}$. In this context, it is very easy to generalize to schemes with vector fields the various constructions  done for $\diffSp A$. 
\medskip

Then, we compare the three different sheaves that have been defined over $\diffSp A$ :
\begin{itemize}
\item  the restricted sheaf $\mathcal{O}_{\diffSp A}^{\text{(Keigher)}}$, defined in \Cite{KeigherQuasiAffine} ;
\item  the sheaf $\mathcal{O}_{\diffSp A}^{\text{(Kovacic)}}$, defined \emph{\`a la Hartshorne} in \Cite{CarraFerro2} and used by Kovacic in several papers ;
\item the sheaf  $\mathcal{O}_{\diffSp A}^{\text{(CF)}}$ defined by Carr\`a Ferro in \Cite{KolchinSchemes}. 
\end{itemize}
We prove that 
\[\mathcal{O}_{\diffSp A}^{\text{(Keigher)}}\cong \mathcal{O}_{\diffSp A}^{\text{(Kovacic)}}\]
for any differential ring $A$. For the Carr\`a Ferro sheaf, we prove (see Theorem \Ref{Comparison.Const.Gal}):

\begin{theoremeintro}
Let $X$ be a reduced $\Q$-scheme with a vector field. Then, the Carr\`a Ferro sheaf and the Keigher sheaf have the same constants:
$$
\fa \, U \text{open in $X^{\vec{\mathscr{V}}}$}, \qquad C \big ( \mathcal{O}_{X^{\vec{\mathscr{V}}}}^{\mathrm{(Keigher)}}(U) \big ) \simeq C \big ( \mathcal{O}_{X^{\vec{\mathscr{V}}}}^{\mathrm{(CF)}}(U) \big ).
$$
\end{theoremeintro}

\smallskip

\noindent The main ingredient of the proof is the following proposition. 


   \begin{propositionintro}
Let $X$ be a reduced $\Q$-scheme with a vector field $\vec{\mathscr{V}}$. Let $U$ be an open set of $X$. Then, for every $f\in C\left ( \mathcal{O}_{X} ( U ) \right )$, there exists a unique $\widetilde{f}$ in $C\left ( \mathcal{O}_{X}(U^\delta) \right )$ such that $\widetilde{f}_{\mid U}=f$. \\[0mm]

\noindent Furthermore, the extension map
\[
\mathrm{ext}_{U\to U^\delta}:\left\{\fonctionb
{C \left ( \mathcal{O}_{X}(U) \right )}
{C \left ( \mathcal{O}_{X}(U^\delta) \right )}
{f}
{\widetilde{f}}
\right.
\]
is an isomorphism of rings, whose inverse ${C \left ( \mathcal{O}_{X}(U^\delta) \right )} \longto {C \left ( \mathcal{O}_{X}(U) \right )}$  is the restriction map. 
\end{propositionintro}
   
\smallskip

\noindent The proof of this proposition relies on the following result of commutative algebra, proved in this paper.

\begin{propositionintro}
Let $A$ be a differential ring and $S$ a multiplicative subset of $A$. For $s\in S$ and $i\in \N$, we have
$$
\left . \begin{array}{r}
\left ( \displaystyle\frac{a}{s} \right )'=0 \text{ in $S^{-1}A$} \\
s^{(i)}\in S
\end{array} \right \} \quad \impl \quad \frac{a}{s}=\frac{a^{(i)}}{s^{(i)}}\quad\text{ in $S^{-1}A$}.
$$
\end{propositionintro}

   \bigskip
   

\bigskip

\section{Vector fields, leaves and trajectories for schemes}

We start this paper with some classical and elementary facts about vector fields on smooth manifolds. This will be a good motivation for our definition for schemes.

\subsection{Vector fields }In the case of smooth manifolds, we can define global vector fields in various ways. Given $M$ such a manifold:
\begin{itemize}
\item[a)] If the tangent bundle $TM$ has already been defined, we can say that a global vector field is a section $s$ of the canonical projection $\pi: TM \longto M$.
\item[b)] It is equivalent to consider a map
$$
\partial: \mathscr{C}^\infty\left ( M, \R \right ) \longto \mathscr{C}^\infty\left ( M, \R \right )
$$
that is $\R$-linear and such that
$$
\fa f,g\in \mathscr{C}^{\infty}\left ( M, \R \right ), \qquad \partial \left ( fg \right )=f\partial(g)+\partial(f)g.
$$
In  other words, global vector fields can also be seen as $\R$-derivations of the $\R$-algebra $\mathscr{C}^\infty\left ( M, \R \right )$. The derivation $\partial_s$ associated to a  section $s$ of $TM \longto M$ is defined by
$$
\partial_{s} (f):=\quad
 \fonctionb{M}{\R}{p}{df_{p}\bullet s(p)}
$$
for all $f\in  \mathscr{C}^\infty\left ( M, \R \right )$. Let us give more details about this construction, for the convenience of the reader. First, $df_p$, by definition, is a linear application from $T_p M$ to $\R$. By $df_{p}\bullet s(p)$, we simply mean the image of the tangent vector $s(p)\in T_p M$ by the map $df_p$. So, $\partial_s(f)$ is a real function defined on $M$. We will not prove that it is actually a smooth function, but let us see now why the map $\partial_s : \mathscr{C}^\infty\left ( M, \R \right ) \longto  \mathscr{C}^\infty\left ( M, \R \right ) $ is a derivation. It is a consequence of the following:
\begin{lemma}
 Let $M$ be a smooth manifold and $f,g\in \mathscr{C}^\infty\left ( M, \R \right )$. Let $p\in M$ and $\vec{v}\in T_p M$. Then, 
 $$
 d_p(fg)\bullet \vec{v}=f(p)\cdot d_p(g)\bullet \vec{v} + g(p)\cdot d_p(f)\bullet \vec{v}.
 $$
\end{lemma}
\noindent This lemma is easily deduced from the same result for $\R^n$, $n\geq 1$.
\item[c)] Actually, given a global vector field, one gets a $\R$-derivation $\partial_{U}$ of $\mathscr{C}^{\infty}\left ( U, \R \right )$ for all open set $U$ of $M$. Moreover, these maps are compatible with the restriction maps. Thus, one can attach to a global vector field a $\R$-derivation of the structure sheaf $\mathscr{O}_{M}$ of $M$. 
\end{itemize}

\noindent This motivates the definition:

\smallskip
\begin{definition}
Let $X$ be a scheme. A \emph{vector field} $\vec{\mathscr{V}}$ on $X$ is a derivation of the structure sheaf $\mathscr{O}_{X}$ of $X$.
\end{definition}

\smallskip

\remarks{For the convenience of the reader, let us give more details. A vector field on $X$ is therefore given by the following data:  
\begin{itemize}
 \item  for each open set $U$ of $X$, a derivation $\partial_U$ of the ring $\mathscr{O}_X(U)$ such that
 \item for all open sets $U,V$ of $X$ such that $U \subset V$, the following diagram commutes
 $$
 \xymatrix@C=1.5cm@R=8mm{
 \mathscr{O}_X(V) \ar[d]_-{\rho_{V,U}} \ar[r]^-{\partial_V} & \mathscr{O}_X(V) \ar[d]^-{\rho_{V,U}}  \\
  \mathscr{O}_X(U)\ar[r]^-{\partial_U} & \mathscr{O}_X(U) 
 }.
 $$
\end{itemize}

\point This definition stands in the more general setting of ringed spaces.

\point A scheme has always at least one vector field: the zero-vector field.

\point
If $(X, \mathscr{O}_{X})$ is a scheme, then it is equivalent to consider a vector field $\vec{\mathscr{V}}$ on $X$ or to endow the sheaf $\mathscr{O}_{X}$ with a structure of sheaf of differential rings: $( X, \mathscr{O}_{X}, \vec{\mathscr{V}} )$ is then what we will call a \emph{differentially ringed space}.

\point In  \Cite{EGA44}, given a $S$-scheme $X$, Grothendieck  defines the \emph{tangent bundle of $X/S$}. It is a $S$-scheme, denoted by $T_{X/S}$, with a $S$-morphism to $X$:
$$
\xymatrix@R=5mm{
T_{X/S}\ar[d]^\pi\\X.
}
$$
He proves that the $S$-section of $\pi$ correspond to the $\mathscr{O}_{S}$-derivations of $\mathscr{O}_{X}$. So, in the case where $X$ is viewed as a $\Z$-scheme, one gets a correspondence between the sections of $\pi : T_{X}\longto X$ and the group of vector fields of $X$. 
The $\mathscr{O}_{X}$-module of $S$-sections of $\pi$ is the dual of $\Omega_{X/S}^1$. We will denote it by $\mathscr{T}_{X/S}$ (or by $\mathscr{T}_{X}$ when $S=\Sp \Z$).}

\bigskip

\noindent 

\bigskip

\subsection{Morphisms and category }If $\mathscr{X}=(X, \vec{\mathscr{V}})$ and $\mathscr{Y}=(Y, \vec{\mathscr{W}})$ are two schemes with vector fields, a morphism $f: \mathscr{X}\longto \mathscr{Y}$ will be a morphism $f: X \longto Y$ of schemes such that, for all open set $U$ of $Y$, the diagram
$$
\xymatrix@C=2cm{
\mathscr{O}_{X}\left ( f^{-1}\left ( U \right ) \right )\ar[d]^{\partial_{\vec{\mathscr{V}}, f^{-1}(U)}} & \mathscr{O}_{Y}\left ( U \right )\ar[l]_-{f^{\#}_{U}} \ar[d]^{\partial_{\vec{\mathscr{W}}, U}} \\
\mathscr{O}_{X}\left ( f^{-1}\left ( U \right ) \right )& \mathscr{O}_{Y}\left ( U \right )\ar[l]_-{f^{\#}_{U}} 
}
$$
commutes. In other words, $f$ is a morphism of schemes that is a morphism of differentially ringed spaces. The category of schemes with vector fields will be denoted by $\dSch$. Intuitively, as it will be seen in Proposition \Ref{Prop:Compa:Morph}, a morphism $f: \mathscr{X}\longto \mathscr{Y}$ pushes the vector field of $\mathscr{X}$ onto the vector field of $\mathscr{Y}$.

\bigskip

\subsection{Schemes with vector fields and differential rings }\label{AffinesDiffSchemes}If $(A, \partial)$ is a differential ring, then the scheme $\Sp A$ can be canonically endowed with a vector field, that will be denoted here by  $\vec{\mathscr{V}}_{A}$. Let's build it.

 First, we remark that it is sufficient to define $\vec{\mathscr{V}}_{A}$ on a basis of open sets. The collection $\left \{D(f), f\in A \right \}$ is a basis of $\Sp A$, where $D(f)$ denotes the open set $D(f):=\left \{ \mathfrak{p}\in \Sp A \tq f\notin \mathfrak{p} \right \}$. The ring of sections on $D(f)$ (for the structure sheaf of $\Sp A$) is the localisation $A_f$. In other words:
$$
\mathscr{O}_{\Sp A} \left ( D(f) \right ) \simeq A_f \qquad \text{naturally}.
$$
So, to define $\vec{\mathscr{V}}_{A}$, we just need to define a ``compatible'' collection of derivations of the rings $A_f$, $f\in A$. Now, we know that the derivation $\partial$ of $A$ induces a unique ``natural'' derivation $\partial_f$ of the localisation $A_f$. We leave to the reader the ``natural'' task to verify that these derivations $\partial_f$ form a ``compatible'' collection. That is how the vector field $\vec{\mathscr{V}}_{A}$ of $\Sp A$ is defined. We will denote this scheme with vector field by $\dSp A$. Actually, one obtains a functor
$$
\dSp: (\dAnn)^\mathrm{op}\longto \dSch.
$$
We could also have defined the schemes with vector field as differentially ringed spaces locally isomorphic to $\dSp A_{i}$'s.

\bigskip

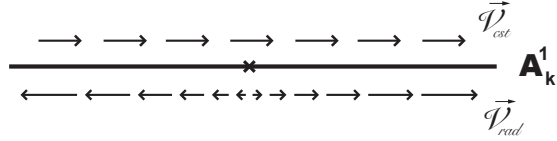
\begin{figure}[h]
\begin{center}
\begin{tikzpicture}
\tikzset{>=angle 90}
	\draw[line width=2pt,] (-4.5,0) -- (4.5,0);
	\draw[fill] (0,0) circle (3pt);
\foreach \x in {-2,...,2}{
\draw[line width=1pt,shift={(\x*1.5,0)},->] (-0.5,.5) -- +(1,0);
}

\foreach \x in {0,1,...,4}{
\pgfmathsetmacro{\u}{(\x+1)*(\x+2)/10-.2+.1}
\pgfmathsetmacro{\l}{0.2*(\x+1)}
\draw[line width=1pt,shift={(\u,0)},->] (0,-0.5) -- +(\l,0);
\draw[line width=1pt,shift={(-\u,0)},->] (0,-0.5) -- +(-\l,0);
}

\node () at (5.2,0) {$\A{1}{k}$};
\node () at (4.2, 0.55) {$\overrightarrow{\mathscr{V}}_{\text{cst}}$};
\node () at (4.6,-0.55) {$\overrightarrow{\mathscr{V}}_{\text{rad}}$};
\end{tikzpicture}
\end{center}
\caption{The vector fields of $\A{1}{k}$ associated to the derivations $\partial^{\text{cst}}$ and $\partial^{\text{rad}}$.}\label{Dessin:examples}
\end{figure}

\example{
Let $k$ be a field and $A=k[x]$. The derivation $\partial^{\text{cst}}$ of $
A$ defined by
${\partial^{\text{cst}}}_{\mid k}=0$ and $\partial^{\text{cst}}(x)=1$
corresponds to the constant vector field of $\A{1}{k}$. The derivation $\partial^{\text{rad}}$ defined by
${\partial^{\text{rad}}}_{\mid k}=0$ and $\partial^{\text{rad}}(x)=x$
corresponds to the radial vector field, as pictured in Figure \Ref{Dessin:examples}.
}

	\bigskip

As in the non-differential case, one has the following proposition, whose proof is left to the reader: 
\smallskip

\begin{proposition}\label{adjonction}
The functors
$$
\boite{\xymatrix@C=2cm{
(\dAnn)^\mathrm{op}\ar@{}|-{\top}[r]\ar@<1.6mm>[r]^-{\dSp }& \ar@<1.6mm>[l]^-{\mathscr{O}\left ( - \right )}\dSch
}}
$$
form an adjunction: $\mathscr{O}\left ( - \right )$ is a left adjoint to $\dSp$.
\end{proposition}

\smallskip

\noindent In particular, the category of affine schemes with vector field is antiequi\-valent to the category of differential rings. This 
allows us to describe the vector fields of $\A{n}{k}$ and $\Pn{n}{k}$, as follows.

\bigskip

\examples{\emph{Vector fields on $\A{n}{k}$.} Let $k$ be a ring. Let $\vec{\mathscr{V}}$ be a \mbox{vector} field defined on $\A{n}{k}$, and constant on $k$. Then, $\vec{\mathscr{V}}$  corresponds to a $k$-derivation of $k[X_{1}, \ldots, X_{n}]$. Such a derivation $\partial$ is fully determined by the elements $\partial X_{1}, \ldots, \partial X_{n}$. Hence, the abelian group of vector fields on $\A{n}{k}$ is isomorphic to $(k[X_{1}, \ldots, X_{n}],+)^n$

\bigskip

\point \emph{Vector fields on $\Pn{n}{k}$.}
Let $k$ be a ring. Then, the vector fields defined on $\Pn{n}{k}$ ($n\geq 1$) and constant on $k$ all come from linear vector fields of $\A{n+1}{k}$. This means, precisely, that for any vector field $\vec{\mathscr{V}}$ defined on $\Pn{n}{k}$ and constant on $k$, there exists a matrix $A\in M_{n+1}(k)$ such that the morphism
$$
\pi : \xymatrix{(\A{n+1}{k}\setminus \left \{0\right \},  \vec{\mathscr{V}}_{A})\ar[r]& (\Pn{n}{k}, \vec{\mathscr{V}})}
$$
is compatible with the vector fields, where $\pi : \A{n+1}{k}\setminus \left \{0\right \}\longto \Pn{n}{k}$ denotes the canonical projection and where $\vec{\mathscr{V}}_{A}$  denotes the linear vector field of $\A{n+1}{k}$ induced by the derivation
$$
\partial_{A} :{k[X_{0}, \ldots, X_{n}]}\longto {k[X_{0}, \ldots, X_{n}]}
$$
defined by 
$$
\partial_{A} \left ( \!\! \begin{array}{c}
X_{0}\\ \vdots \\\ X_{n}
\end{array} \! \!\right ) := A\left (\! \! \begin{array}{c}
X_{0}\\ \vdots \\\ X_{n}
\end{array} \! \!\right ) .
$$

\bigskip

Indeed, the Euler exact sequence can be written
$$
0 \longto \mathscr{O}_{\Pn{n}{k}}\longto \mathscr{O}_{\Pn{n}{k}}(1)^{n+1}\longto \mathscr{T}_{\Pn{n}{k}/k }\longto 0
$$
as in Example 8.20.1 of \Cite{hartshorne}. Hence, one gets an exact sequence in cohomology:
$$
H^0\!\big ( \Pn{n}{k}, \mathscr{O}_{\Pn{n}{k}}(1)^{n+1} \big )\longto
H^0\!\big ( \Pn{n}{k}, \mathscr{T}_{\Pn{n}{k}/k} \big )\longto
H^1\!\big ( \Pn{n}{k}, \mathscr{O}_{\Pn{n}{k}} \big ).
$$
But, one knows that $H^1\!\big ( \Pn{n}{k}, \mathscr{O}_{\Pn{n}{k}} \big )=0$ (see for instance \Cite{AlgGeoArithCurves}). So, the map $H^0\!\big ( \Pn{n}{k}, \mathscr{O}_{\Pn{n}{k}}(1)^{n+1} \big )\longto
H^0\!\big ( \Pn{n}{k}, \mathscr{T}_{\Pn{n}{k}/k} \big )$ is surjective. Let us write down expli\-citly what is this map. To a family $(L_{0}, \ldots, L_{n})$ of linear forms in $X_{0}, \ldots, X_{n}$, it associates the vector field of $\Pn{n}{k}$, defined on each standard open set $U_{i}= \Sp k\left [ X_{0}/X_{i}, \ldots,{X_{n}}/{X_{i}}\right ]$ by:
$$
\partial (X_{k}/X_{i} ) = \frac{ L_{k}\cdot X_{i}- L_{i} \cdot X_{k}}{{X_{i}}^2}.
$$
So, given a vector field $\vec{\mathscr{V}}$ of $\Pn{n}{k}$, one obtains the required matrix $A$ by considering the coefficients of the linear forms $L_{0}, \ldots, L_{n}$.
}

\bigskip

\subsection{Tangent vectors associated to vector fields }Let $X$ be a scheme. Now, we are going to explain how to associate to a vector field $\vec{\mathscr{V}}$ on $X$ and to an element $x\in X$ a (Zariski)  tangent vector $\vec{\mathscr{V}}(x)\in T_{x}X$. First, it is easy to check that a vector field $\vec{\mathscr{V}}$, \ie a derivation $\partial$ of $\mathscr{O}_{X}$, induces a derivation $\partial_{x}$ of the local ring $\mathscr{O}_{X,x}$. We denote by $\mathfrak{M}_{x}$, as usual, the maximal ideal of $\mathscr{O}_{X,x}$, and $\kappa(x):= \mathscr{O}_{X,x}/\mathfrak{M}_{x}$. We then consider the linear map
$$
\fonctionb
{\mathfrak{M}_{x}}
{\kappa(x)}
{f}
{(\partial_{x}f)(x)}.
$$
This map sends elements of $\mathfrak{M}_{x}^2$ to zero, since
$$
\partial_{x}(fg)(x)=\left ( ( \partial_{x}f  )g + f (\partial_{x} g)\right )(x)=0
$$
for $f, g \in \mathfrak{M}_{x}$. Hence,  we get a map
$$
\mathfrak{M}_{x}/ \mathfrak{M}_{x}^2\longto \kappa(x),
$$
which is $\kappa(x)$-linear: in other words, we get a element of $T_{x}X$ the Zariski tangent space of $X$ in $x$. We denote this element by $\vec{\mathscr{V}}(x)$. We then have the formula:

\smallskip

\begin{proposition}\label{Prop:Compa:Morph}
Let $\mathscr{X}=(X, \vec{\mathscr{V}})$ and $\mathscr{Y}=(Y, \vec{\mathscr{W}})$ be two schemes with vector fields. Let $f: \mathscr{X}\longto \mathscr{Y}$ be a morphism. Then
$$
\fa x\in X, \qquad T_{x} f \bullet \vec{\mathscr{V}}(x)=i_{x}\circ\vec{\mathscr{W}}(f(x)).
$$
where $i_{x} : \kappa(f(x))\longto \kappa(x)$ is the inclusion of residual fields induced by $f$.
\end{proposition}

\smallskip

\pproof{Since the definition of $\vec{\mathscr{V}}(x)$ is local, it is sufficient to prove this statement when $X$ and $Y$ are affine. So, let $(A, \partial_{A})$ and $(B, \partial_{B})$ be two differential rings, and let $\varphi: (A, \partial_{A})\longto (B, \partial_{B})$ be a morphism. We denote by $f: \dSp B \longto \dSp A$ the corresponding morphism of schemes with vector fields. Let $x\in \Sp B$, \ie let $\mathfrak{p}_{x}$ be a prime ideal of $B$. The image of $x$ by $f$ is $\mathfrak{p}_{y}:=\varphi^{-1}(\mathfrak{p}_{x})$. The morphism $\varphi$ induces, by localisation, an arrow
$$
\widehat{\varphi}: A_{\mathfrak{p}_{y}}\longto B_{\mathfrak{p}_{x}}.
$$
Better, if we denote by $\mathfrak{M}_{x}$ and $\mathfrak{M}_{y}$ the maximal ideals of $A_{\mathfrak{p}_{x}}$ and $A_{\mathfrak{p}_{y}}$, $\varphi$ induces morphisms
$$
\overline{\varphi}: \mathfrak{M}_{y}/ \mathfrak{M}_{y}^2 \longto\mathfrak{M}_{x}/ \mathfrak{M}_{x}^2 \qquad 
\aand 
\qquad i_{x}: B_{\mathfrak{p}_{y}}/\mathfrak{M}_{y}  \longto A_{\mathfrak{p}_{x}}/\mathfrak{M}_{x},
$$
this latter being injective.
Now, the tangent vector $\vec{\mathscr{V}}(x)$ corresponds to the morphism
$$
\partial_{x}: \fonctionb
{\mathfrak{M}_{x}/ \mathfrak{M}_{x}^2}
{B_{\mathfrak{p}_{x}}/ \mathfrak{M}_{x}}
{\psi}
{\partial_{B} (\psi)\text{ mod. }\mathfrak{M}_{x}},
$$
and there is a similar description of $\vec{\mathscr{W}}(y)$.
The image of $\vec{\mathscr{V}}(x)$ by the differential $T_{x}f$ is the map $\partial_{y}$ making the diagram
$$
\xymatrix{
\ar[dr]_-{\partial_{y}}\mathfrak{M}_{y}/\mathfrak{M}_{y}^2\ar[r]^-{\overline{\varphi}}& \mathfrak{M}_{x}/\mathfrak{M}_{x}^2 \ar[d]^-{\partial_{x}} \\
&A_{\mathfrak{p}_{x}}/\mathfrak{M}_{x}
}
$$
commute.
Let $\psi\in \mathfrak{M}_{y}$. We have: 
\begin{align*}
\partial_{y}(\psi \text{ mod. }\mathfrak{M}_{y}^2) &=\partial_{x} \left ( \overline{\varphi}(\psi \text{ mod. }\mathfrak{M}_{y}^2) \right ) = \partial_{x} \left ( \widehat{\varphi}(\psi) \text{ mod. }\mathfrak{M}_{x}^2 \right ) \\
&=\partial_{B}\left ( \widehat{\varphi}(\psi) \right ) \text{ mod. }\mathfrak{M}_{x} \\
&= \widehat{\varphi} \left ( \partial_{A} (\psi) \right )\text{ mod. }\mathfrak{M}_{x} \\
&= i_{x} \left ( \partial_{A}(\psi) \text{ mod. }\mathfrak{M}_{y} \right ).
\end{align*}
In other words, $T_{x}f \bullet \vec{\mathscr{V}}(x)=i_{x}\circ \vec{\mathscr{W}}(f(x))$. }

\bigskip

\subsection{Leaves }We are now able to define leaves: 

\smallskip

\begin{definition}
Let $\mathscr{X}=(X, \vec{\mathscr{V}})$ be a scheme with a vector field. Let $\eta \in X$. We say that \emph{$\eta$ is a leaf of $\mathscr{X}$} (or \emph{a leaf for $\vec{\mathscr{V}}$}) when $\vec{\mathscr{V}}(\eta)=0$. The set of leaves of $\mathscr{X}$ will be denoted by $X^{\vec{\mathscr{V}}}$.
\end{definition}

\smallskip

Let us check that the leaves of $\dSp A$ correspond to the differential prime ideals of $A$, when $A$ is a differential ring. Let $\mathfrak{p}$ be a prime ideal of $A$. Let's assume that $\mathfrak{p}$ is a leaf of $\dSp A$. Let $f\in \mathfrak{p}$: from $\vec{\mathscr{V}}(\mathfrak{p})=0$, we deduce that the image of $f$ under the map
$$
\mathfrak{p}A_{\mathfrak{p}}=\mathfrak{M}_{\mathfrak{p}}\longto \mathfrak{M}_{\mathfrak{p}}/\mathfrak{M}_{\mathfrak{p}}^2\stackrel{\partial_{A}} {\longto}
\mathfrak{M}_{\mathfrak{p}}/\mathfrak{M}_{\mathfrak{p}}^2
\longto A_{\mathfrak{p}}/\mathfrak{M}_{\mathfrak{p}}
$$
is zero. Hence, $\partial_{A}(f)\in \mathfrak{M}_{\mathfrak{p}}=\mathfrak{p}A_{\mathfrak{p}}$ and so, $f\in \mathfrak{p}$: the ideal $\mathfrak{p}$ is differential. Conversely, one can check that if $\mathfrak{p}$ is a differential ideal, then it is a leaf of $\dSp A$. This fundamental remark shows that
$$
X^{\vec{\mathscr{V}}}\subset X
$$
is the exact non-affine analogue of
$$
\diffSp A \subset \Sp A.
$$

\bigskip

\examples{The scheme  $\A{1}{k}$ endowed with the constant vector field has only one leaf: its generic point $\eta$. With the radial vector field, it has two leaves: the closed point $0$ and $\eta$.

\bigskip

\point Let's consider the ring $A=k[X_{1}, \ldots, X_{n}]$ with a $k$-derivation $\partial$. The derivation $\partial$ is characterized par the elements
$$
P_{1}:=\partial(X_{1}) \qquad P_{2}:=\partial (X_{2}) \qquad \cdots \qquad P_{n}:=\partial(X_{n}).
$$
One can check that the corresponding vector field $\vec{\mathscr{V}}$ satisfies, for all $x_{1}, \ldots, x_{n}\in k$: 
$$
\vec{\mathscr{V}} \left ( x_{1}, \ldots x_{n} \right ) = \left ( \! \begin{array}{c}
P_{1}(x_{1}, \ldots, x_{n}) \\
\vdots \\
P_{n}(x_{1}, \ldots, x_{n}) \\
\end{array} \! \right ).
$$
Let's take $n=2$ (we denote $A=k[x,y]$) with the derivation
$$
\partial(x)=-2y \qquad \aand \qquad \partial(y)=3x^2.
$$
By a simple computation, one checks that the prime ideal $\eta_{c}=(x^3+y^2-c)$ is differential, for all $c\in k$: consequently, $(\eta_{c})_{c\in k}$ is a family of leaves.

\bigskip

\point As noticed by Buium in \Cite{buium} (Lemma (2.1) of Chapter 1), if $X$ is a $\Q$-scheme and if $F$ is an irreducible closed set of $X$, then the generic point $\eta_{F}$ of $F$ is always a leaf, for  any vector field $\vec{\mathscr{V}}$. 

\bigskip

\point Let $k$ be a ring. Any vector field on $\Pn{n}{k}$ vanishes on a closed point --- a kind of analogue of the \emph{hairy ball theorem}. Indeed, as explained in \RefPar{AffinesDiffSchemes}, if $\vec{\mathscr{V}}$ is a vector field of $\Pn{n}{k}$ constant on $k$, then there exists $A\in M_{n+1}(k)$ such that $\vec{\mathscr{V}}$ comes from the vector field of $\A{n+1}{k}$ defined by
$$
\partial\left ( \!\! \begin{array}{c}
X_{0}\\ \vdots \\\ X_{n}
\end{array} \! \!\right ) := A\left (\! \! \begin{array}{c}
X_{0}\\ \vdots \\\ X_{n}
\end{array} \! \!\right ) .
$$
Now, let $K$ be a residual field of $k$, \ie let $\varphi : k\longto K$ be a surjective morphism. There exists a finite extension of fields $K \longto L$ such that the matrix $A$, when viewed in $L$, has an eigenvector $\vec{v}$. This  implies that the image of $\vec{v}$ under the map $\pi_{L} : \A{n+1}{L}\longto \Pn{n}{L}$, denoted by $x_{L}$, and which is clearly a closed point, is a leaf for $\vec{\mathscr{V}}_{L}$ --- in other words, $\vec{\mathscr{V}}_{L}$ vanishes on $x_{L}$. By Proposition \Ref{Prop:Compa:Morph}, one knows that the image of $x_{L}$ under the map
$$
f : \Pn{n}{L}\longto \Pn{n}{K}\longto \Pn{n}{k},
$$
denoted by $x_{k}$, will also be a leaf for $\vec{\mathscr{V}}$. Hence, we just need to see why $x_{k}$ is a closed point. This comes from the following facts:
\begin{itemize}
\item[---] First, since $L/K$ is finite, $\Sp L \longto \Sp K$ is a proper map and so is $\Pn{n}{L}\longto \Pn{n}{K}$. In  particular, it is a closed map.
\item[---] Second, since $\Sp K \longto \Sp k$ is a closed immersion, the morphism $\Pn{n}{K}\longto \Pn{n}{k}$ is also a closed immersion. In particular, it is a closed map.
\item[---] So, $\Pn{n}{L}\longto \Pn{n}{k}$ is a closed map, and sends closed points to closed points: $x_{k}$ is a closed point.
\end{itemize}
}

\bigskip


\subsection{Trajectory of a point }Now, let $\mathscr{X}=(X, \vec{\mathscr{V}})$ be in $\dSch$. We would like to associate to any $x\in X$ ``its algebraic trajectory under the vector field $\vec{\mathscr{V}}$''. This is possible, thanks to the following theorem, which is an analogue for schemes of the Cauchy-Peano theorem: 

\smallskip

\begin{theorem}\label{Theo:CauchyLip:schemes}
Let $\mathscr{X}=(X, \vec{\mathscr{V}})$ be a scheme with a vector field, defined over $\Q$. Let $x\in X$. Then, the ordered set
$$
\left \{\eta \in X \, \middle | \,\begin{array}{c}
\eta\leadsto  x \\
\eta \text{ is a leaf of }\mathscr{X}
\end{array}\right \}
$$
has a least element. We denote this element by $\text{Traj}_{\vec{\mathscr{V}}}(x)$ and call it the trajectory of $x$ (under $\vec{\mathscr{V}}$).
\end{theorem}

\smallskip

\remark{Here, the order that we consider is: $z\geq y$ if and only if $y\in \overline{\left \{z\right \}}$. In this case, we say that \emph{$z$ is a generization of $y$} and that \emph{$y$ is a specialization of $z$}; we denote $z\leadsto y$. The properties of this order are classical (see \Cite{EGA1}, Chapter 0, \textbf{(2.1.1)}). For instance, open sets are stable under generization and, dually, closed sets are stable under specialization.}

\bigskip

\pproof{We keep the notations of the theorem. Let $x\in X$ and $U$ an affine neighborhood of $x$. Since all the  generizations of $x$ are elements of $U$, one can assume $X$ affine. So, let $A$ be a differential $\Q$-algebra and $\mathfrak{p}$ a prime ideal of $A$. Since the generization order is the opposite of the inclusion order on ideals, one needs to prove that
$$
\left \{ \mathfrak{q}\in \Sp A \, \middle | \, \begin{array}{c}
\mathfrak{q}\text{ is a \emph{differential} prime ideal} \\
\mathfrak{q}\subset \mathfrak{p}
\end{array}\right \}
$$
has a greatest element. Let's consider, as in \Cite{KeigherSpecialRings}, the set
$$
\mathfrak{p}_{\#}:=\left \{f\in A \tq \fa n\geq 0, \quad f^{(n)}\in \mathfrak{p}\right \}.
$$
Keigher's Proposition 1.5 says that $\mathfrak{p}_{\#}$ is a prime ideal (it's there that $\Q\subset A$ is needed). It is then easy to check that $\mathfrak{p}_{\#}$ is the required ideal. We will see further a proof of the primality of $\mathfrak{p}_{\#}$ is a more general context. Let's remark that, in any case, $I_{\#}$ is the greatest differential ideal contained in $I$.
}

\bigskip

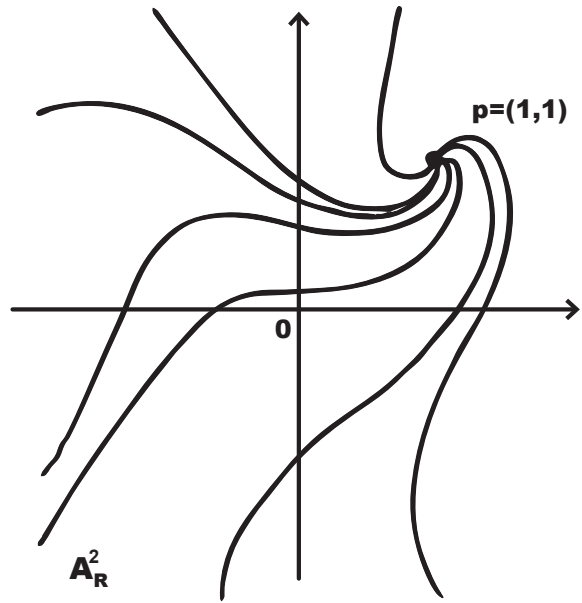
\begin{figure}[h]
\begin{center}
\begin{tikzpicture}[x=2cm, y=2cm]
\tikzset{>=angle 90}
\newcommand{\xmin}{-2}
\newcommand{\ymin}{-2}
\newcommand{\xmax}{2}
\newcommand{\ymax}{2}
\draw[thick,->,color=black] (\xmin,0) -- (\xmax,0) node[below right] {};

\draw[thick,->,color=black] (0,\ymin) -- (0,\ymax) node[right] {};


\clip (\xmin,\ymin) rectangle (\xmax,\ymax);
\foreach \depx/\depy in {-2/-1.5,1/-2,-2/1,-1/2,1/2,-2/-1,-.5/-2}{
\pgfmathsetmacro\x{\depx}
\pgfmathsetmacro\y{\depy}
\pgfmathsetmacro\pre{0.01} 
	\pgfplotsforeachungrouped  \n in {1,2,...,500}{
			\pgfmathsetmacro\nx{\x+\pre*(1-\x*\y*\y)}
			\pgfmathsetmacro\ny{\y+\pre*(\x*\x-\y*\y*\y)}
			\draw[thick] (\x,\y) -- (\nx,\ny) ;
			\pgfmathsetmacro\x{\nx}
			\pgfmathsetmacro\y{\ny}
	}

}

\node () at (-1.5,-1.8) {$\A{2}{\R}$};
\draw[fill,color=black] (1,1) circle (2pt) node[yshift=7mm, xshift=5mm] { $p=(1, 1)$};
\end{tikzpicture}
\end{center}
\caption{Smooth leaves of the vector field defined by $
\partial x=1-xy ^2 \aand{} \partial y=x^2-y^3$.}\label{Dessin:Jouanolou}
\end{figure}

\examples{For all leaf $\eta \in X,\,\, \text{Traj}_{\vec{\mathscr{V}}}(\eta)=\eta$.

\smallskip

\point Since  $X=\A{1}{k}$ endowed with the constant vector field has only $\eta$ as a leaf, one has:
$\fa x\in X, \, \text{Traj}_{\vec{\mathscr{V}}}(x)=\eta.$
If we consider the radial vector field, 
the trajectory of all $x\in X$ but $0$ is $\eta$.

\smallskip

\point Let's consider the vector field on $\A{2}{\C}$ defined by
$$
\partial x=1-xy ^2 \qquad \aand \qquad \partial y=x^2-y^3,
$$
whose smooth real leaves are drawn in Picture \Ref{Dessin:Jouanolou}. Jouanolou proved in \Cite{jouanolou} that no non-constant smooth leaf of this vector field is algebraic. Thus, the leaves for this vector field are just $\eta$ the generic point and the point $(1,1)$.}

\subsection{Properties of the trajectory} First, we prove that the map $\text{Traj}_{\vec{\mathscr{V}}}$ is ``compatible'' with morphisms of $\dSch$, namely:

\smallskip

\begin{proposition}
Let $\mathscr{X}=(X, \vec{\mathscr{V}})$ and $\mathscr{Y}=(Y, \vec{\mathscr{W}})$ be two $\Q$-schemes with vector fields, and let $f: \mathscr{X}\longto \mathscr{Y}$ be a morphism. Then, for all $x\in X$,
$$
f \big (\text{Traj}_{\vec{\mathscr{V}}}(x)\big )= \text{Traj}_{\vec{\mathscr{W}}}(f(x)).
$$
\end{proposition}

\smallskip

\pproof{
By considering affine neighborhoods of $f(x)$ and $x$, it suffices to prove this proposition in  the affine case. Hence, let $A$ and $B$ be two $\Q$-differential algebras, let $\varphi: A\longto B$ be a morphism of differential rings. Let $\mathfrak{p}$ be a prime ideal of $B$. We want to prove that
$$
\varphi^{-1} \big ( \mathfrak{p}_{\#} \big )=\big ( \varphi^{-1}\mathfrak{p} \big )_{\#}.
$$
But, 
\begin{align*}
\varphi^{-1}\left (  \mathfrak{p}_{\#} \right )  &= \left \{x\in A \tq \varphi(x)\in \mathfrak{p}_{\#}\right \} = \left \{x\in A \tq \fa n\in \N,\, {\varphi(x)}^{(n)}\in \mathfrak{p}\right \} \\
&= \left \{x\in A \tq \fa n\in \N,\, {\varphi\big ( x^{(n)} \big )}\in \mathfrak{p}\right \} = \left \{x\in A \tq \fa n\in \N,\, x^{(n)} \in \varphi^{-1} \mathfrak{p}\right \} \\
&= \left ( \varphi^{-1} \mathfrak{p} \right )_{\#},
\end{align*}
which concludes the proof.
}

\bigskip

The trajectory defines a map
$$
\text{Traj}_{\vec{\mathscr{V}}}: X \longto X^{\vec{\mathscr{V}}}.
$$
Since $X^{\vec{\mathscr{V}}} \subset X$, it is possible to endow the set $X^{\vec{\mathscr{V}}}$ of leaves with the topology induced by the Zariski topology. Then:

\bigskip

\begin{proposition}\label{Traj.cont.open}
Let  $X$ be $\Q$-scheme endowed with a vector field $\vec{\mathscr{V}}$. Then, 
$\text{Traj}_{\vec{\mathscr{V}}}: X \longto X^{\vec{\mathscr{V}}}$
 is continuous and open.
\end{proposition}

\smallskip

\pproof{First, let's show that it is continuous. Since this property is local, let's assume that $X=\dSp A$, with $A$ a differential ring. Let $U=V \cap X^{\vec{\mathscr{V}}}$ be an open set of $X^{\vec{\mathscr{V}}}$, where $V$ is a Zariski open set of $X$. Let $I$ be an ideal of $A$ such that
$V=X\setminus V(I)$.
Let's prove that
$$
\left ( \text{Traj}_{\vec{\mathscr{V}}} \right )^{-1} U =X\setminus V(\langle I \rangle),
$$
where $\langle I \rangle$ denotes the differential ideal generated by $I$. Let $\mathfrak{p}$ be a prime ideal of $A$. Then, one has
\begin{align*}
\text{Traj}_{\vec{\mathscr{V}}}(\mathfrak{p})\in U & \ssi \mathfrak{p}_{\#}\in U \\
& \ssi \mathfrak{p}_{\#}\in V \\
& \ssi I \subset \mathfrak{p}_{\#}.
\end{align*}
But, the latter is equivalent to $\langle I \rangle \subset  \mathfrak{p}$.  Indeed, if $I\subset \mathfrak{p}_{\#}$, since $\mathfrak{p}_{\#}$ is a differential ideal, one has $\langle I \rangle \subset  \mathfrak{p}_{\#}$ and since $\mathfrak{p}_{\#}\subset \mathfrak{p}$, one has indeed  $\langle I \rangle \subset  \mathfrak{p}$. On the other hand, if  $\langle I \rangle \subset  \mathfrak{p}$, since $\mathfrak{p}_{\#}$ is greatest differential ideal contained in $\mathfrak{p}$, one has $\langle I \rangle \subset  \mathfrak{p}_{\#}$ and so $I \subset   \mathfrak{p}_{\#}$. This proves indeed that $
\left ( \text{Traj}_{\vec{\mathscr{V}}} \right )^{-1} U =X\setminus V(\langle I \rangle)$.

\bigskip

Let's prove now that $\text{Traj}_{\vec{\mathscr{V}}}$ is open. Let $X$ be a $\Q$-scheme endowed with a vector field $\vec{\mathscr{V}}$. Let $U$ be an open set of $X$. Since, for all $\eta\in X^{\vec{\mathscr{V}}}$, one has $\Traj{V}(\eta)=\eta$, it is easy to check that
$$
\Traj{V}(U)=U \cap X^{\vec{\mathscr{V}}}.
$$
Hence, the map $\Traj{V}$ is indeed open.
}

\subsection{The case when the schemes are no more defined over~$\Q$}

A crucial hypothesis in Theorem \Ref{Theo:CauchyLip:schemes} is that the schemes have to be defined over $\Q$. This comes from the fact that, when differentiating $f^n$, one gets $n\cdot f^{n-1}f'$: if $n$ can be simplified, much more things can be done. In general, Theorem \Ref{Theo:CauchyLip:schemes} is false when the schemes are not defined over $\Q$. Nevertheless, it is possible to solve this problem by defining \emph{Hasse-Schmidt} vector fields. Recall that, when $A$ is a ring, a \emph{Hasse-Schmidt derivation of $A$} is a family $D=(D_{i})_{i\geq 0}$ of map $A\longto A$ satisfying 
\begin{itemize}
\item[\quad (i)] for all $i\geq 0$, $D_{i}: A\longto A$ is an additive map, and $D_{0}=\Id_{A}$.
\item[\quad (ii)] \emph{the generalized Leibniz rule}: for all $i$ and all $f,g\in A$: 
$$
D_{i}(fg)=\sum_{k+\ell=i}D_{k}(f)D_{\ell}(g).
$$
\item[\quad (iii)] \emph{iterativity}: for all $i,j\geq 0$,
$D_{i}\circ D_{j}=\binom{i+j}{i}D_{i+j}$.
\end{itemize}
If $A$ is a $\Q$-algebra, then, there is a one-to-one correspondance between derivations of $A$ and Hasse-Schmidt derivations of $A$ given by
$$
\xymatrix@C=1.5cm{
\partial \ar@{|->}[r]& D:=\left (\displaystyle{ \frac{\partial^i}{i !}}\right )_{i\geq 0}.
}
$$
Subsequently, if $X$ is a scheme, we call \emph{Hasse-Schmidt vector field of $X$} any Hasse-Schmidt derivation of the structure sheaf $\mathscr{O}_{X}$, \ie any family $(D_{U})_{U}$ of compatible Hasse-Schmidt derivations $\mathscr{O}_{X}\left ( U \right )\longto \mathscr{O}_{X}\left ( U \right )$ for all open set $U$. Let's now define what would be a leaf for a Hasse-Schmidt vector field. The situation is more complicated than for classical vector fields. For any Hasse-Schmidt derivation $\mathscr{D}$  of $\mathscr{O}_{X}$ and any $x\in X$, it is possible to consider the restriction ${\mathscr{D}}_{x}$ of $\mathscr{D}$ to the local ring $\mathscr{O}_{X,x}$: it is a Hasse-Schmidt derivation of $\mathscr{O}_{X,x}$, and we denote $\mathscr{D}_{x}=(\mathscr{D}_{x,i})_{i\geq 0}$. For all $i\geq 1$,  the map
$$
\text{ev}_{x,i}: \fonctionb{\mathfrak{M}_{x}^{i}/\mathfrak{M}_{x}^{i+1}}{\kappa(x)}{f}{\mathscr{D}_{x,i}(f)(x)}
$$
is well defined, since
$$
\fa f\in \mathfrak{M}_{x}^{i+1},\qquad \mathscr{D}_{x,i}(f)(x)=0.
$$
Indeed, if $A$ is a ring and if $D=(D_{0}, D_{1}, \ldots)$ is a Hasse-Schmidt derivation of $A$, the generalized Leibniz rule generalises to
$$
\fa i\geq 0, \, \fa p\geq 1, \qquad 
D_{i}\left ( f_{1}f_{2}\cdots f_{p} \right )
=
\sum_{\begin{subarray}{c}
\ell_{1}, \ldots, \ell_{p}\geq 0 \\[0.5mm]
\ell_{1}+\cdots +\ell_{p}=i
\end{subarray}} D_{\ell_{1}}(f_{1})\cdots D_{\ell_{p}}(f_{p}),
$$
for all $f_{1}, \ldots, f_{p}\in A$.
We say that $x$ is a leaf for this Hasse-Schmidt derivation if the maps $\text{ev}_{x,i}$ are zero for all $i\geq 1$. Then we have:

\smallskip

\begin{theorem}
Let $X$ be a scheme endowed with a Hasse-Schmidt vector field $\vec{\mathscr{V}}$. Let $x\in X$. Then, the ordered set
$$
\left \{\eta \in X \, \middle | \,\begin{array}{c}
\eta\leadsto  x \\
\eta \text{ is a leaf for }\vec{\mathscr{V}}
\end{array}\right \}
$$
has a least element.
\end{theorem}

\smallskip

The proof of this theorem is based on the following proposition:

\smallskip

\begin{proposition}
Let $A$ be a ring and $D=(D_{i})_{i\geq 0}$ a Hasse-Schmidt derivation of $A$. Let $\mathfrak{p}$ be a prime ideal of $A$. Then,
$$
\mathfrak{p}_{\#}:=\left \{f\in A \tq \fa i\geq 0,\, D_{i}(f)\in \mathfrak{p}\right \}
$$
is a prime ideal invariant by $D$.
\end{proposition}

\smallskip

\pproof{
The set $\mathfrak{p}_{\#}$ is clearly stable under addition. If $f\in \mathfrak{p}_{\#}$ and $\lambda\in A$, then the generalized Leibniz rule proves that $\lambda f \in \mathfrak{p}_{\#}$. Furthermore, the iterativity of $D$ proves that for all $i\geq 0$, the ideal $\mathfrak{p}_{\#}$ is stable under $D_{i}$. Let's prove that $\mathfrak{p}_{\#}$ is a prime ideal. Let $f,g\in A$ such that $f,g\notin \mathfrak{p}_{\#}$. Thus, let $i_{0}$ and $j_{0}\geq 0$ be the least integers such that
$$
D_{i_{0}}(f)\notin \mathfrak{p} \qquad \aand \qquad D_{j_{0}}(g)\notin \mathfrak{p}.
$$
Let's prove that $fg\notin  \mathfrak{p}_{\#}$ by considering
$$
D_{i_{0}+j_{0}}(fg)=\sum_{k+\ell=i_{0}+j_{0}}D_{k}(f)D_{\ell}(g).
$$
In this sum, the terms split in three parts: the $D_{k}(f)D_{\ell}(g)$'s for $k<i_{0}$, which are in $\mathfrak{p}$ by definition of $i_{0}$, the $D_{k}(f)D_{\ell}(g)$'s for $\ell<j_{0}$, which are in $\mathfrak{p}$ for the same reason, and finally $D_{i_{0}}(f)D_{j_{0}}(g)$. This latter isn't in  $\mathfrak{p}$ for $\mathfrak{p}$ is a prime ideal. Thus, $D_{i_{0}+j_{0}}(fg) \notin \mathfrak{p}$ and so, $fg\notin \mathfrak{p}_{\#}$. This proves that $\mathfrak{p}_{\#}$ is a prime ideal.
}

\bigskip

\section{The Carr\`a Ferro topology, the Carr\`a Ferro sheaf and the Keigher sheaf}

In this section, we reinterpret the paper \Cite{KolchinSchemes}, with the help of vector fields, leaves and trajectories. This new approach allow us to generalize the constructions of Carr\`a Ferro to the non-affine case and, much more important, to get a geometric understanding of these latter.

\subsection{Invariant closed and open sets }We begin this section by defining  invariant closed and open sets. For the sake of simplicity, we stick to $\Q$-schemes with vector fields but all what follows should work for schemes with Hasse-Schmidt derivations.

\smallskip

\begin{definition}
Let $X$ be a $\Q$-scheme with a vector field $\vec{\mathscr{V}}$. A closed set $F$ of $X$ will be said \emph{invariant under $\vec{\mathscr{V}}$} when
$$
\fa x\in F,\quad  \text{Traj}_{\vec{\mathscr{V}}}(x)\in F.
$$
An open set $U$ of $X$ will be said \emph{invariant under $\vec{\mathscr{V}}$} when the closed set $X\setminus U$ is.
\end{definition}

We now prove an analogue of Theorem \Ref{Theo:CauchyLip:schemes} for open sets: 

\smallskip

\begin{proposition}
Let $X$ be a $\Q$-scheme endowed with a vector field $\vec{\mathscr{V}}$. Let $U$ be an open set of $X$. Then, the set 
$$
\left \{
V \,\middle |\,
\begin{array}{c}
U\subset V \\
V \text{ is an invariant open set of }{X}
\end{array}
\right \}
$$
has a least element. We denote it by $U^\delta$, and call it the \emph{invariant open set associated to $U$}.
\end{proposition}

\remark{Of course, dually, there also exists a greatest invariant closed set included in $F$, when $F$ is a closed set of $X$.}

\medskip

\pproof{We keep the notations of the statement. Let's consider the set 
$$V_{0}=\left \{x\in X \tq \Traj{V}(x)\in U \right \}.$$
Since the map $\Traj{V}$ is continuous, $V_{0}$ is an open set of $X$. Furthermore, $U \subset V_{0}$, for we always have $\Traj{V}(x) \leadsto x$, and for open sets are stable under generization. Now, let's check that $V_{0}$ is invariant: let $x\notin V_{0}$. Then, $\Traj{V}(x)\notin V_{0}$, since $\Traj{V} \left ( \Traj{V}(x) \right )=\Traj{V}(x)\notin  U$. Last, let's prove that $V_{0}$ is the least such set. Let $V\supset U$ be an invariant open set and $x\in V_{0}$. If $x\notin V$, then, by invariance, one would have $\Traj{V}(x)\notin V$. But, by definition of $x\in V_{0}$, one has $\Traj{V}(x)\in U$ and thus 
$\Traj{V}(x)\in V$. This is absurd. Hence, $x\in V$ and so, $V_{0}\subset V$.
}

\bigskip

\noindent If $A$ is a differential ring, if $X=\dSp A$ and if $U$ is the open set of $X$ defined by an ideal $I$, then $U^\delta$ is the open set defind by the differential ideal $\left \langle I \right \rangle$. Indeed, one has
$$
U^\delta = \left ( \text{Traj}_{\vec{\mathscr{V}}} \right )^{-1} U =X\setminus V(\langle I \rangle),
$$
as it has been shown in the proof of Proposition \Ref{Traj.cont.open}.

\subsection{The Carr\`a Ferro topology of $\mathscr{X}$} We now prove that the invariant open sets form a topology: 

\smallskip

\begin{proposition}
Let $\mathscr{X}=(X, \vec{\mathscr{V}})$ be a $\Q$-scheme endowed with a vector field. Let $(U_{i})_{i\in I}$ be a family of open sets of $X$. Then: 
$$\left ( \bigcup_{i\in I} U_{i} \right )^\delta = \bigcup_{i\in I} {U_{i}}^\delta \quad \text{and, when $I$ is finite,}\quad
 \left ( \bigcap_{i\in I} U_{i} \right )^\delta = \bigcap_{i\in I} {U_{i}}^\delta .
$$
In particular, the invariant open sets of $X$ form a topology of $X$. We call it \emph{the Carr\`a Ferro topology fo $\mathscr{X}$.}
\end{proposition}

\smallskip

\pproof{This comes from the fact that for any map $f : E \to F$, $f^{-1}$ commutes with unions and intersections, applied to $f=\Traj{V}$.
}

\smallskip

Consequently, the subset ${X}^{\vec{\mathscr{V}}}$ of $X$ can be endowed with two induced topologies: the one induced by Zariski, and the one induced by Carr\`a Ferro. They are the same: 

\smallskip

\begin{proposition}
Let $\mathscr{X}=(X, \vec{\mathscr{V}})$ be $\Q$-scheme endowed with a vector field. Then, the Zariski topology of $X$ and the Carr\`a Ferro topology of $\mathscr{X}$ induce the same topology on $X^{\vec{\mathscr{V}}}$.
\end{proposition}

\smallskip

\pproof{Since, the Carr\`a Ferro topology is a subtopology of the Zariski topo\-logy, it suffices to prove that, if $U$ is Zariski open set of $X$, then, there exists an invariant open set $V$ of $X$ such that:
$$
U \cap X^{\vec{\mathscr{V}}} = V \cap X^{\vec{\mathscr{V}}}.
$$
It suffices to take $V:=U^\delta$.
}

\subsection{The Carr\`a Ferro sheaf  and the Keigher sheaf on $\mathscr{X}$} Now, we would like to equip the topological space $X^{\vec{\mathscr{V}}}$ with a sheaf. For this, we have three possibilities. 
\begin{itemize}
\item[a)] First, if we denote by $X_{\text{Zar}}$ the scheme $X$ \emph{endowed with the Zariski topology}, then the inclusion map
$$
i_{\text{Zar}}: X^{\vec{\mathscr{V}}}\longto X_{\text{Zar}}
$$
is a continuous map. Since $X_{\text{Zar}}$ comes with the scheme-structure scheaf $\mathscr{O}_{X}$, one can consider the pull-back of $\mathscr{O}_{X}$ by $i_{\text{Zar}}$. In other words, one can consider the restriction of $\mathscr{O}_{X}$ to the subspace $X^{\vec{\mathscr{V}}}$. It is a sheaf denoted by 
$$
(i_{\text{Zar}})^{-1} \mathscr{O}_{X}
$$
and defined as the sheaf associated to the preasheaf
$$
U \mapsto \limind \limits_{
\begin{subarray}{c}
 V\text{ open in $X$} \\
 \text{and }U \subset V
\end{subarray}
}\, \mathscr{O}_{X}\left ( V \right ).
$$
Indeed, this latter is not always a sheaf. This sheaf is naturally a sheaf of differential $\Q$-algebras.

\smallskip

\item[b)] Second, we can do the same but with the Carr\`a Ferro topology instead of the Zariski one. So, if we denote by $X_{\text{CF}}$ the scheme $X$ equipped \emph{with the Carr\`a Ferro topology}, it is still possible to consider the inclusion map
$$
i_{\text{CF}}: X^{\vec{\mathscr{V}}}\longto X_{\text{CF}}: 
$$
it is also a continuous map. The sheaf $\mathscr{O}_{X}$, defined on $X_{\text{Zar}}$, induces naturally a sheaf on $X_{\text{CF}}$, which we still denote by $\mathscr{O}_{X}$. Thus, similarly, one can consider the sheaf
$$
(i_{\text{CF}})^{-1} \mathscr{O}_{X}.
$$

\smallskip

\item[c)] Third, there is another sheaf that one can define on $X^{\vec{\mathscr{V}}}$. Indeed, since $\Traj{V}: X_{\text{Zar}} \longto X^{\vec{\mathscr{V}}}$ is a continuous map and since $X_{\text{Zar}}$ comes with the sheaf $\mathscr{O}_{X}$, one can consider the push-forward of $\mathscr{O}_{X}$ by $\Traj{V}$. It is a sheaf denoted by
$$
(\Traj{V})_{*}\mathscr{O}_{X}
$$
and whose definition is simplier than for the pull-back: if $U$ is a open set of $X^{\vec{\mathscr{V}}}$, one has, by definition
$$
(\Traj{V})_{*}\mathscr{O}_{X} \left ( U \right ):= \mathscr{O}_{X}\left (  (\Traj{V})^{-1} U\right ).
$$
\end{itemize}

\begin{notation}
When $U$ is a open set of $X^{\vec{\mathscr{V}}}$, we denote
$$
U_{\Delta}:=  (\Traj{V})^{-1} U = \left \{x\in X \tq \Traj{V}(x)\in U\right \}.
$$
\end{notation}

\noindent It is easy to check that $U_{\Delta}$ is an invariant open set of $X$. With this notation, we have
$(\Traj{V})_{*}\mathscr{O}_{X} \left ( U \right )  =  \mathscr{O}_{X}\left (   U_{\Delta} \right )$. We have: 

\bigskip

\begin{proposition}
Let $\mathscr{X}=(X, \vec{\mathscr{V}})$ be $\Q$-scheme endowed with a vector field. Then, 
$$
(i_{\text{CF}})^{-1} \mathscr{O}_{X}=(\Traj{V})_{*}\mathscr{O}_{X}.
$$
\end{proposition}

\bigskip

\pproof{We keep the notations of the proposition. Let $U$ be an open  set of $X^{\vec{\mathscr{V}}}$. We will prove that $\mathscr{O}_{X}\left ( U_{\Delta} \right )$ is an inductive limit of the $\mathscr{O}_{X}\left ( V \right )$, for $V$ invariant open set of $X$ such that $U \subset V$. So, let $V$ be an invariant open set containing $U$. Then, $U_{\Delta} \subset V$. 
Hence, the restrictions form  a bunch of compatible maps
$$
 \psi_{V}: \mathscr{O}_{X}\left ( V \right )\longto \mathscr{O}_{X}\left ( U_{\Delta} \right ).
$$
Let's prove that these maps make $\mathscr{O}_{X}\left ( U_{\Delta} \right )$ an inductive limit. It's easy. Let $A$ be a differential ring, equipped with compatible maps $\varphi_{V}: \mathscr{O}_{X}\left ( V \right )\longto A$ for all invariant open set $V$ containing $U$. In  particular, there is a map 
$$
f:=\varphi_{U_{\Delta}}: \mathscr{O}_{X}\left ( U_{\Delta} \right ) \longto A.
$$
What we want to prove is that, for every $V$, the diagram
$$
\xymatrix{
\ar@/_6mm/[rr]_{\varphi_{V}}\mathscr{O}_{X}\left ( V \right )\ar[r]^{\psi_{V}}& \mathscr{O}_{X}\left ( U_{\Delta} \right )\ar[r]^f& A
}
$$
commutes. This follows from the compatibility of the family $(\varphi_{V})_{V}$.
}

\bigskip

\begin{definition}
Let $\mathscr{X}=(X, \vec{\mathscr{V}})$ be $\Q$-scheme endowed with a vector field. The \emph{Keigher sheaf} on $X^{\vec{\mathscr{V}}}$ is
$$
\mathscr{O}_{X^{\vec{\mathscr{V}}}}^{\mathrm{(Keigher)}}:= (i_{\text{Zar}})^{-1} \mathscr{O}_{X}.
$$ 
The \emph{Carr\`a Ferro sheaf} on $X^{\vec{\mathscr{V}}}$ is
$$
\mathscr{O}_{X^{\vec{\mathscr{V}}}}^{\mathrm{(CF)}}:= (i_{\text{CF}})^{-1} \mathscr{O}_{X}=(\Traj{V})_{*}\mathscr{O}_{X}.
$$ 
\end{definition}

\bigskip

With these definitions, Corollary 2.4 of \Cite{KolchinSchemes} generalizes to the following: 

\smallskip

\begin{proposition}
Let $\mathscr{X}=(X, \vec{\mathscr{V}})$ be $\Q$-scheme endowed with a vector field. Then,
$$
\Gamma \big ( X^{\vec{\mathscr{V}}}, \mathscr{O}_{X^{\vec{\mathscr{V}}}}^{\mathrm{(CF)}}  \big ) = \Gamma(X, \mathscr{O}_{X}).
$$
In particular, if $A$ is $\Q$-differential algebra, 
$$
\Gamma \big ( \diffSp A, \mathscr{O}_{\diffSp A}^{\mathrm{(CF)}}  \big )
=A.
$$
\end{proposition}

\smallskip

\pproof{By definition, 
$$
\Gamma \big ( X^{\vec{\mathscr{V}}}, \mathscr{O}_{X^{\vec{\mathscr{V}}}}^{\mathrm{(CF)}}  \big ) = \Gamma((X^{\vec{\mathscr{V}}})^\delta, \mathscr{O}_{X}).
$$
But, it is clear that $(X^{\vec{\mathscr{V}}})_{\Delta}=X$ and thus, the result follows.
}

\smallskip

\subsection{The Kovacic sheaf }When $X$ is affine, a third sheaf has been studied. Although it has been defined for the first time by Carr Ferro in \Cite{CarraFerro2}, we call it the \emph{Kovacic sheaf}. Indeed, in a series of papers \Cite{KovacicDiffSchemes, KovacicGlobalSections, KovacicGaloisDiffSchemes1, KovacicGaloisDiffSchemes2}, Kovacic intensively uses and studies this sheaf. Here is its definition:

\smallskip

\begin{definition}
Let $A$ be a differential ring and $U$ an open set of $\diffSp A$. The \emph{Kovacic sheaf} $\mathscr{O}^{\textrm{(Kov)}}_{\diffSp A}$, is defined by
\begin{gather*}
\mathscr{O}^{\textrm{(Kov)}}_{\diffSp A} \left ( U \right )
:= \\ \\
\left \{s : U \longto \coprod_{\mathfrak{p}\in U} A_{\mathfrak{p}} \middle | \begin{array}{lc}
(i) &\fa \mathfrak{p}\in U,\,  s(\mathfrak{p})\in A_{\mathfrak{p}} \\ \\
&
 \xt \left ( U_{i} \right )_{i\in I}\text{ open covering of }U, \\
(ii)  &\xt (a_{i})_{i\in I}, (b_{i})_{i\in I}\in A^I, \\
 & \fa \mathfrak{p}\in U, \, \fa i \in I, \quad \mathfrak{p}\in U_{i}\impl (b_{i}\notin \mathfrak{p}\aand s(\mathfrak{p})=\frac{a_{i}}{b_{i}})
\end{array}\right \}.
\end{gather*}
\end{definition}

We will prove further that $\mathscr{O}^{\textrm{(Kov)}}_{\diffSp A}$ and 
$\mathscr{O}^{\textrm{(Keigher)}}_{\diffSp A}$ are isomorphic.
\bigskip

\section{Extension of constants}

In this section, we prove the following result, which will be our main tool to compare the Keigher sheaf and the Carr\`a Ferro sheaf. If $A$ is a differential ring, we denote by $C(A)$ the ring of constants of $A$.

\begin{proposition}\label{proposition:extension}
Let $\mathscr{X}=(X, \vec{\mathscr{V}})$ be a reduced $\Q$-scheme endowed with a vector field. Let $U$ be an open set of $X$. Then, for every $f\in C\left ( \mathscr{O}_{X} ( U ) \right )$, there exists a unique $\widetilde{f}$ in $C\left ( \mathscr{O}_{X}(U^\delta) \right )$ such that $\widetilde{f}_{\mid U}=f$. \\[0mm]

\noindent Furthermore, this extension map
$$
\mathrm{ext}_{U\to U^\delta}:\fonctionb
{C \left ( \mathscr{O}_{X}(U) \right )}
{C \left ( \mathscr{O}_{X}(U^\delta) \right )}
{f}
{\widetilde{f}}
$$
is an isomorphism of rings, whose inverse ${C \left ( \mathscr{O}_{X}(U^\delta) \right )} \longto {C \left ( \mathscr{O}_{X}(U) \right )}$  is the restriction map. 
\end{proposition}

\bigskip

\subsection{Constants in localized rings }In order to prove Proposition \Ref{proposition:extension}, we need to study the properties of constant elements in differential rings of the form $S^{-1}A$. If $x=a/s$ is such an element, by differentiating $x$, one gets
$$
\frac{a's-s'a}{s^2}=0\qquad \text{in $S^{-1}A$}.
$$
One would like to derive from this, identities such as
$$
a's-s'a=0 \qquad \text{and so} \qquad \frac{a}{s}=\frac{a'}{s'}\qquad \text{and so} \qquad \fa i\in \N, \quad \frac{a}{s}=\frac{a^{(i)}}{s^{(i)}}.
$$
Unfortunately, these latters are false, since we don't have $a's-s'a=0$ but only $\xt t\in S, \,\,t\cdot (a's-s'a)=0$, and since the elements $s^{(i)}$ do not necessarily belong to $S$. Nevertheless, when $s^{(i)}\in S$, we do have $a/s=a^{(i)}/s^{(i)}$ in $S^{-1}A$. This is what tells us the following proposition.

\bigskip

\begin{proposition}\label{Constantes.Dans.Un.Anneau.Fractions}
Let $A$ be a differential ring. Let $\theta, a,b\in A$ such that
$$
\theta \cdot \left ( a'b-ab' \right )=0.
$$
1) Then, for all $N\in \N_{\geq 3}$ and for all $0\leq i \leq N$, one has
\begin{gather*}
\theta\cdot \left ( a'b-ab' \right )=0 \\
\theta^2 \cdot \left ( a''b-ab'' \right )=0 \\
b^{N-1}\theta  ^{N}\cdot \left ( b^{(i)}a^{(N-i)}-a^{(i)}b^{(N-i)} \right )=0.
\end{gather*} 
2) In particular, when $S$ is a multiplicative subset of $A$, if $s\in S$ and if $i\in \N$, one has
$$
\left . \begin{array}{r}
\left ( \frac{a}{s} \right )'=0 \text{ in $S^{-1}A$} \\
s^{(i)}\in S
\end{array} \right \} \quad \impl \quad \frac{a}{s}=\frac{a^{(i)}}{s^{(i)}}\quad\text{ in $S^{-1}A$}.
$$
\end{proposition}

\bigskip

In order to prove this proposition, we need the following lemma :

\bigskip

\begin{lemma}\label{Lemme.Derivee.Anneau.Frac}
Let $A$ be a differential ring and let $t, A_{1}, A_{2}, B_{1}, B_{2}\in A$. Then, 
\begin{gather*}
t \cdot \left ( A_{1}B_{2}-B_{1}A_{2} \right )=0 \\
\Downarrow \\
t^2\cdot \left ( {B_{2}}^2\cdot \left ( B_{1}A_{1}'-A_{1}B_{1}' \right )-{B_{1}}^2\cdot \left ( A_{2}'B_{2}-B_{2}'A_{2} \right ) \right )=0.
\end{gather*}
\end{lemma}

\bigskip

\pproofbis{of Lemma \Ref{Lemme.Derivee.Anneau.Frac}}
{Let $A$ be a differential ring and $t, A_{1}, A_{2}, B_{1}, B_{2}\in A$. Let's denote $\Theta=t \cdot \left ( A_{1}B_{2}-B_{1}A_{2} \right )$. A simple computation shows that
\begin{gather*}
t B_{1}B_{2}  \frac{\partial \Theta}{\partial x}  -   t B_{1}B_{2}'\cdot  \Theta - t B_{1}'B_{2} \Theta -t'B_{1}B_{2}\Theta \\
= \\
t^2 \cdot \left ( {B_{2}}^2\cdot \left ( A_{1}'B_{1}-A_{1}B_{1}' \right ) -{B_{1}}^2\cdot \left ( A_{2}'B_{2}-A_{2}B_{2}' \right )\right ).
\end{gather*}
Hence, when $\Theta=0$, one gets the required identity.
}

\bigskip

Now, we can prove Proposition \Ref{Constantes.Dans.Un.Anneau.Fractions} :

\bigskip

\pproofbis{of Proposition \Ref{Constantes.Dans.Un.Anneau.Fractions}}{We keep the notations of the proposition. In particular, we assume that
$\theta \cdot \left ( a'b-ab' \right )=0$. We denote
$$
E_{N,i} \quad := \quad b^{N-1}\theta  ^{N}  \cdot \left ( b^{(i)}a^{(N-i)}-a^{(i)}b^{(N-i)} \right ).
$$
Let's begin by showing the assertion $\emph{1)}$: we want to prove\begin{align*}
&\theta\cdot \left ( a'b-ab' \right )=0 \\
&\theta^2 \cdot \left ( a''b-ab'' \right )=0 \\
\fa N\geq 3,\quad  \fa 0\leq i \leq N, \qquad &b^{N-1}\theta  ^{N}\cdot \left ( b^{(i)}a^{(N-i)}-a^{(i)}b^{(N-i)} \right )=0.
\end{align*} The first identity is our assumption;  one gets the second one by differentiating the first one and by multiplying it by $\theta$. For the bunch of next identities, we proceed by induction. For $N=3$, let's remark that, when differentiating the second identity and multiplying it by $\theta$, one gets:
\begin{equation}
\label{Deux.Eq.Derivee}
\theta^3\cdot \left (\left ( a''b'-a'b'' \right ) + \left ( a'''b-ab''' \right ) \right )=0.
\end{equation}
But, by applying Lemma \Ref{Lemme.Derivee.Anneau.Frac} with $t=\theta$, $A_{1}=a$, $B_{2}=b'$, $B_{1}=b$ et $A_{2}=a'$, one gets 
$$
b^2 \theta^2 \cdot \left ( b' a'' -a'b'' \right )=0;
$$
hence, in particular, one has
$$
b^2 \theta^3 \cdot  \left ( b' a'' -a'b'' \right )=0 \qquad \text{ and, with (\ref{Deux.Eq.Derivee}),}\qquad b^2 \theta^3 \cdot \left (  a'''b-ab''' \right )=0
$$
Now, let's assume the assertion \emph{1)} true for $n\leq N$ and let's show it for  $N+1$. First, a simple computation shows that
$$
b \theta  \cdot \frac{\partial  E_{N,i}}{\partial x} =E_{N+1, i+1}+E_{N+1,i}.
$$
Thus, for $0\leq i \leq N$, one has $E_{N+1, i+1}+E_{N+1,i}=0$. A consequence of these identities is that, if there exists $i_{0}$ such that  $E_{N+1, i_{0}}=0$ then all the $E_{N+1, i}$ are zero. Indeed, in that case, one would have
$$
\underbrace{\begin{blockarray}{cccccccccc}
      \begin{block}{(cccccccccc)}
        1& 1 & 0  &\cdots \\
        0&1&1       \\
         \vdots            \\
                &&& \ddots&\ddots     \\
                     \\
                     \\
                                          &&&\cdots &&0&1&1&0\\
                     &&&& &&0&1&1\\
                     && \cdots&0&1&0&\cdots&&\\
        \end{block} 
        &&&&i_{0}
        \end{blockarray}
}_{A}        \cdot \left (\! \!\begin{array}{c}
E_{N+1,0} \\
E_{N+1,1} \\
\vdots \\
E_{N+1,N+1} \\
\end{array}\! \!\right )=0.
$$
But, developing along the last line, one finds that $\det A=\left ( -1 \right )^{N+i_{0}}$ and thus that $A$ is invertible. So, it suffices to find $i_{0}$ such that $E_{N+1, i_{0}}=0$. We consider two cases.
If $N+1=2k$ is  even, then one has
\begin{align*}
E_{N+1, k}&=b^{N}\theta^{N+1} \cdot \left (  b^{(k)}a^{(\left ( N+1 \right )-k)}-a^{(k)}b^{(\left ( N+1 \right )-k)} \right ) \\
&=b^{N}\theta^{N+1}\cdot \left (  b^{(k)}a^{(k)}-a^{(k)}b^{k)} \right )=0,
\end{align*}
and we can conclude. If $N+1=2k+1$ is odd, we know, by the induction assumption, that
$$
E_{k,0}=b^{k-1}\theta^{k}\cdot \left ( a^{(k)}b+b^{(k)}a \right )=0.
$$
Then, if we use Lemma \Ref{Lemme.Derivee.Anneau.Frac} with the data 
$$
t=\left ( b^{k-1}\theta^{k} \right ) \quad A_{1}=a^{(k)} \quad B_{2}=b \quad B_{1}=b^{(k)}\quad A_{2}=A,
$$
we get
$$
\left ( b^{k-1}\theta^{k}  \right )^{2}\cdot \left ( b^2\cdot \left ( b^{(k)}a^{(k+1)}-a^{(k)}b^{(k+1)} \right )+{b^{(k)}}^2\cdot \left ( a'b-b'a \right ) \right )=0.
$$
So, given $\theta\cdot \left ( a'b-b'a \right )=0$, we get
$$
b^{2k}\theta^{2k}\cdot \left ( b^{(k)}a^{(k+1)}-a^{(k)}b^{(k+1)} \right )=0.$$
Mulitplying it by $\theta$, we get $E_{N+1, k}=0$ --- and so, all the $E_{N+1, i}$ are zero.

\bigskip

Now, let's move to the assertion \emph{2)}.
It is an easy consequence of \emph{1)}. Indeed, let $S$ be a multiplicative subset of $A$ and let $(a,s)\in A\times S$ such that 
$$
\left ( \frac{a}{s} \right )'=0 \qquad \text{ in $S^{-1}A$.}
$$
It means that there exists $\theta\in S$ such that
$
\theta\cdot \left ( a's-s'a \right )=0.
$
Let's assume now that  $i\in \N$ verifies $s^{(i)}\in S$. The identity $E_{i,0}=0$ that we have just shown tells us that
$$
\underbrace{\left ( s^{i-1}\theta^i \right )}_{\in S}\cdot \left ( a^{(i)}s -as^{(i)}\right )=0.
$$
For $s^{(i)}\in S$, this implies 
$$
\frac{a}{s}=\frac{a^{(i)}}{s^{(i)}} \qquad \text{ in $S^{-1}A$.}
$$}

\subsection{A lemma on stalks and trajectories }We will also need the following:

\begin{lemma}\label{Nulle.Dans.Traj.Nilpo.Dans.Point}
Let $( X, \vec{\mathscr{V}} )$ be a $\Q$-scheme endowed with a vector field. Let $x\in X$, let $U$ be an open neighborhood of  $x$ and let $f\in \mathscr{O}_{X}\left ( U \right )$. Then,
\begin{itemize}
\item[{(i)}] \qquad
$f_{\text{Traj}_{\vec{\mathscr{V}}}(x)}=0 \quad \impl \quad \xt \,n\in \N \tq \left ( f_{x} \right )^n=0.$
\item[{(ii)}] \qquad
$(f_{\text{Traj}_{\vec{\mathscr{V}}}(x)}=0 \aand f'=0) \quad \impl \quad f_{x}=0.$
\end{itemize}

\end{lemma}

\bigskip

\remark{This result is false out of the differential context: if $X$ is a scheme, if $x\in X$ and if $\eta\leadsto x$ is a generization of  $x$, then 
$$
f_{\eta}=0\quad \Longrightarrow\hspace*{-3.5ex}/\hspace*{1ex} \quad\xt n\in \N \tq \left ( f_{x} \right )^n =0.
$$
To see this, it suffices to consider the closed subscheme of $\A{2}{\C}$, union of the axes $x=0$ and $y=0$ : $X=\Sp \C[x,y]/(xy)$. In this scheme, the function $y$ is zero in  $\mathscr{O}_{X, \eta_{x}}$ --- where $\eta_{x}$ stands for the generic point of the axe $y=0$ ---
although $y$ is not nilpotent in $\mathscr{O}_{X, (0,0)}$.}

\bigskip

\pproof{For the property we want to show is local, it suffices to prove it for affine schemes. So, let $A$ be a differential ring. To begin with, let's prove the small following result : 
$$
\fa (\theta, f)\in A^2, \qquad \theta f=0 \qquad \impl \qquad (\fa n\in \N, \quad \theta^{(n)}f^{n+1}=0).
$$
We proceed by induction: if $\theta^{(n)}f^{n+1}=0$, by differentiating this identity, one gets
$$
\theta^{(n+1)}f^{n+1} + (n+1)\theta^{(n)}f'f^{n}=0.
$$
By multiplying the latter by $f$, one gets $\theta^{(n+1)}f^{n+2}=0$.
Now let's move to assertion \emph{(i)}: let $\mathfrak{p}$ be a prime ideal of $A$ and let $f\in A$ such that 
$$
f=0 \qquad \text{in $A_{\mathfrak{p}_{\#}}$}.
$$
This means that there exists $\theta\notin \mathfrak{p}_{\#}$ such that $\theta f=0$. But, $\theta\notin \mathfrak{p}_{\#}$ means that there exists $n\in \N$ such that $\theta^{(n)}\notin \mathfrak{p}$. Since, we know that $\theta^{(n)}f^{n+1}=0$, we have
$$
f^{n+1}=0 \qquad \text{in $A_{\mathfrak{p}}$}.
$$
Lastly, let's prove \emph{(ii)}. With the previous notations, we assume, in addition that $f'=0$. From $\theta f=0$, one gets, by induction, that $\theta^{(m)}f=0$ for all $m$. In particular, one has that $\theta^{(n)}f=0$ and so $f=0$ in $A_{\mathfrak{p}}$.
}

\bigskip

\subsection{Proof of Proposition \Ref{proposition:extension}}Now, we come to the proof of our result on the extension of constant sections of the structure sheaf. So, let $X$ be reduced $\Q$-scheme, equipped with a vector field $\vec{\mathscr{V}}$. Let $U$ be an open set of $X$ and let $f\in C \big ( \mathscr{O}_{X}\left ( U \right ) \big )$. We start by proving the unicity of a extension of $f$ to $U^\delta$. So, let 
$\widetilde{f}^{1},  \widetilde{f}^{2}\in C \left ( \mathscr{O}_{X} ( U^\delta ) \right )$ such that  ${\widetilde{f}^{1}}_{\mid U}={\widetilde{f}^{2}}_{\mid U}=f$. Let $x\in U^\delta$. This means that  $\text{Traj}_{\vec{\mathscr{V}}}(x)\in U$. Let's denote $y:=\text{Traj}_{\vec{\mathscr{V}}}(x)$. Thus, one has
$$
{{\widetilde{f}^{1}}}_{y}={{\widetilde{f}^{2}}}_{y}.
$$
Consequently, by lemma \Ref{Nulle.Dans.Traj.Nilpo.Dans.Point}.\emph{(ii)}, one has ${{\widetilde{f}^{1}}}_{x}={{\widetilde{f}^{2}}}_{x}$. Hence,  ${{\widetilde{f}^{1}}}={{\widetilde{f}^{2}}}$, what we wanted to show.

\bigskip

Let's prove now existence of such a extension. Let's assume that we had shown it in the affine case and let's show it in the general case. Let $(\Omega_i)_{i\in I}$ be a basis of open affine sets of $X$. We denote $U_{i}=U \cap \Omega_{i}$ and $f_{i}=f_{\mid U_{i}}$. According to the affine case, one hence has
$$
\widetilde{f_{i}}\in C\left ( \mathscr{O}_{\Omega_{i}}  \big ( {U_{i}}^\delta_{(\subset \Omega_i)} \big )\right )
$$
such that $f_{i}=\widetilde{f_{i}}_{\mid U_{i}}$, where ${U_{i}}^\delta_{(\subset \Omega_i)}$ stands for the invariant open set {of $\Omega_i$} associated to $U_i$:
$$
{U_{i}}^\delta_{(\subset \Omega_i)}:=\left \{ x\in \Omega_i\tq \Traj{V}(x)\in U_i \right \}.
$$
Let's prove that the $\widetilde{f_i}$ patch together, so that one can derive from them  a function $\widetilde{f}$ extending $f$ on $U^\delta$. First,  remark that
$$
\bigcup_{i\in I} {U_{i}}^\delta_{(\subset \Omega_i)} = U^\delta.
$$
This follows from
$$
{U_{i}}^\delta_{(\subset \Omega_i)} = U^\delta \cap \Omega_i.
$$
Now, if $i$ and $j$ are such that $\Omega_j \subset \Omega_i$, since 
${U_{j}}^\delta_{(\subset \Omega_j)} = {U_{i}}^\delta_{(\subset \Omega_i)}  \cap \Omega_j$ and by unicity of the extension, one has
$$
\widetilde{f_j}=(\widetilde{f_{i}})_{\mid \Omega_{j}}.
$$
Finally, if we denote by $\widetilde{f}$ the patching of the $\widetilde{f_i}$, it is clear that $\widetilde{f}'=0$ and that $\widetilde{f}_{\mid U}=f$.

\bigskip

Last, but not least, let's prove the result for affine schemes. For this sake, let's assume that this following lemma is true. We will prove it after.
\smallskip

\begin{lemma}\label{Lemma:Kov}
Let $A$ be a differential reduced ring and let $U$ be an open subset of $\diffSp A$. Let $s\in \mathscr{O}_{\diffSp A}^{\mathrm{(Kov)}}(U)$. Then,
\begin{itemize}
\item[(i)] there exist a Zariski open set $W$ of $\Sp A$, containing $U$,  and $t\in \mathscr{O}_{\Sp A}(W)$ such that for all $x\in U$, the stalk $t_{x}$ equals $s(x)$. Moreover, when $W\subset U_{\Delta}$, this extension $t$ is unique.
\item[(ii)] If, moreover, $s'=0$, this $W$ can be taken to be equal to $U_{\Delta}$: there exists a unique $t\in C \big (\mathscr{O}_{\Sp A}(U_{\Delta}) \big )$ such that
$$
\fa x\in U, \qquad t_x=s(x).
$$

\end{itemize}\end{lemma}

\smallskip

So, let $A$ be a reduced $\Q$-differential algebra, let $V$ be an open set of $X:=\Sp A$, and let $f\in \mathscr{O}_{X}\left ( V \right )$ a section satisfying $f'=0$. We denote $U:= V \cap \diffSp A$. Then, we have $V^\delta=U_{\Delta}$. If we consider the Hartshorne-like \Cite{hartshorne} definition of $f$, then it is clear that $f$ induces on $U$ a constant section $s$ of the Kovacic sheaf. Applying Lemma \Ref{Lemma:Kov}.\emph{(ii)} to $s$, one gets a constant section $\widetilde{f}\in C\left ( \mathscr{O}_{X}(V^\delta) \right )$. We just know that $\widetilde{f}$ and $f$ coincide (in a stalkwise sense) on $U$. But, since $X$ is reduced, by Lemma \Ref{Nulle.Dans.Traj.Nilpo.Dans.Point} this is sufficient to prove that they coincide stalkwisely in $U_{\Delta}$ and so that they are equal. Now, to conclude the proof of Proposition \Ref{proposition:extension}, all that remains is to prove Lemma \Ref{Lemma:Kov}.

\bigskip

\pproofbis{of Lemma \Ref{Lemma:Kov}}
{We keep the notations of the lemma. We start by proving the point \emph{(ii)}. Hence, let $s$ be a constant section of the Kovacic sheaf. It comes with a covering $(U_{i})_{i\in I}$ of $U$ and two families $(a_{i})_{i}$ and $(b_{i})_{i}$ fulfilling the required conditions. The unicity is a consequence of Lemma \Ref{Nulle.Dans.Traj.Nilpo.Dans.Point}.\emph{(ii)}, as for Proposition \Ref{proposition:extension}. For the existence, we use the Hartshorne \Cite{hartshorne} definition of the structure sheaf of $\Sp A$. Hence, we look for 
\begin{itemize}
\item[\emph{\textbf{a)}}] a family $(t(\mathfrak{p}))_{\mathfrak{p}\in U^\delta}$, where $t(\mathfrak{p})\in A_{\mathfrak{p}}$ for each $\mathfrak{p}$, such that
$$
\fa \mathfrak{p}\in U, \quad  t(\mathfrak{p})=s(\mathfrak{p})
\qquad \aand \qquad 
\fa \mathfrak{p}\in U^\delta,\quad t(\mathfrak{p})'=0 .
$$
\item[\emph{\textbf{b)}}] a covering $(\Omega_{\ell})_\ell$ of $U^\delta$ and two families $\left ( \alpha_{\ell} \right )_\ell $ and $\left ( \beta_{\ell} \right )_\ell $ such that
$$
\fa \mathfrak{p}\in U^\delta, \quad \fa \ell, \qquad \left ( \mathfrak{p}\in \Omega_{\ell} \quad\impl \quad  \beta_{\ell}\notin \mathfrak{p} \quad \aand \quad t(\mathfrak{p})=\frac{\alpha_{\ell}}{\beta_{\ell}}\quad \text{ in }A_{\mathfrak{p}} \right ).
$$
\end{itemize} 

\smallskip

\noindent  For the item \emph{\textbf{a)}}, we proceed as follows. Let $\mathfrak{p}\in U^\delta$: this means that $\mathfrak{p}_{\#} \in U$. Hence, one can consider $s(\mathfrak{p}_{\#})$ and  denote
$$
s(\mathfrak{p}_{\#})=\frac{{a_{\mathfrak{p}}}}{{b_{\mathfrak{p}}}},
$$
with ${b_{\mathfrak{p}}}\notin \mathfrak{p}_{\#}$. This means that there exists a integer $n$ such that ${b_{\mathfrak{p}}}^{(n)}\notin \mathfrak{p}$. We will denote the least of these integers by $n_{\mathfrak{p}}$. Then, we define
$$
t(\mathfrak{p}):= \frac{{a_{\mathfrak{p}}}^{(n_{\mathfrak{p}})}}{{b_{\mathfrak{p}}}^{(n_{\mathfrak{p}})}}\in A_{\mathfrak{p}}.
$$
Let's check that this family fulfill the two required conditions. If $\mathfrak{p}\in U$, then $\mathfrak{p}_{\#}=\mathfrak{p}$. So, we want to check that
$$
\frac{{a_{\mathfrak{p}}}}{{b_{\mathfrak{p}}}}= \frac{{a_{\mathfrak{p}}}^{(n_{\mathfrak{p}})}}{{b_{\mathfrak{p}}}^{(n_{\mathfrak{p}})}}\qquad \text{in $A_{\mathfrak{p}_{\#}}$},
$$
but this follows from the point \emph{2)} of Proposition \Ref{Constantes.Dans.Un.Anneau.Fractions}. Let $\mathfrak{p}\in U^\delta$. Since, $s$ is a constant section of the Kovacic sheaf, we know that $s(\mathfrak{p}_{\#})'=0$. That means that
$$
\frac{{a_{\mathfrak{p}}}'{b_{\mathfrak{p}}}-{a_{\mathfrak{p}}}{b_{\mathfrak{p}}}'}{{b_{\mathfrak{p}}}^2}=0\qquad \text{in $A_{\mathfrak{p}_{\#}}$};
$$
thus, there exists $\theta \notin \mathfrak{p}_{\#}$ such that
$$
\theta \cdot ({a_{\mathfrak{p}}}'{b_{\mathfrak{p}}}-{a_{\mathfrak{p}}}{b_{\mathfrak{p}}}')=0.
$$
Now, by point \emph{1)} of Proposition \Ref{Constantes.Dans.Un.Anneau.Fractions}, we know:
\begin{equation}\label{equation.proof}
{b_{\mathfrak{p}}}^{2n_{\mathfrak{p}}}\theta^{2n_{\mathfrak{p}}+1}\cdot
\left ( {a_{\mathfrak{p}}}^{(n_{\mathfrak{p}}+1)}{b_{\mathfrak{p}}}^{(n_{\mathfrak{p}})}-{b_{\mathfrak{p}}}^{(n_{\mathfrak{p}}+1)}{a_{\mathfrak{p}}}^{(n_{\mathfrak{p}})} \right )=0.
\end{equation}
Let's denote $c:={b_{\mathfrak{p}}}^{2n_{\mathfrak{p}}}\theta^{2n_{\mathfrak{p}}+1}$. It is an element that doesn't belong to $\mathfrak{p}_{\#}$ and so, let $m\in \N$ such that $c^{(m)}\notin \mathfrak{p}$. As in the proof of Lemma \Ref{Nulle.Dans.Traj.Nilpo.Dans.Point}, one deduces from (\ref{equation.proof}) that
\begin{align*}
c^{(m)}\cdot \left ( {a_{\mathfrak{p}}}^{(n_{\mathfrak{p}}+1)}{b_{\mathfrak{p}}}^{(n_{\mathfrak{p}})}-{b_{\mathfrak{p}}}^{(n_{\mathfrak{p}}+1)}{a_{\mathfrak{p}}}^{(n_{\mathfrak{p}})} \right )^{m+1}&=0 \\
 \text{and so} \qquad
\left (c^{(m)}\cdot  \left (  {a_{\mathfrak{p}}}^{(n_{\mathfrak{p}}+1)}{b_{\mathfrak{p}}}^{(n_{\mathfrak{p}})}-{b_{\mathfrak{p}}}^{(n_{\mathfrak{p}}+1)}{a_{\mathfrak{p}}}^{(n_{\mathfrak{p}})} \right )\right )^{m+1}&=0.
\end{align*}
Since $A$ is reduced, one has 
$$
c^{(m)}\cdot  \left (  {a_{\mathfrak{p}}}^{(n_{\mathfrak{p}}+1)}{b_{\mathfrak{p}}}^{(n_{\mathfrak{p}})}-{b_{\mathfrak{p}}}^{(n_{\mathfrak{p}}+1)}{a_{\mathfrak{p}}}^{(n_{\mathfrak{p}})} \right )=0.
$$
This implies
$$
\left ( \frac{{a_{\mathfrak{p}}}^{(n_{\mathfrak{p}})}}{{b_{\mathfrak{p}}}^{(n_{\mathfrak{p}})}} \right )'=0\qquad \text{in $A_{\mathfrak{p}}$}.
$$
In other words, it means that $t(\mathfrak{p})'=0$ for all $\mathfrak{p}\in U$.

\bigskip

Let's move now to the item \textbf{\emph{b)}}. For the covering of $U^\delta$, we choose
$$
V_{i,n}:={U_{i}}^\delta \cap D\big ( {b_{i}}^{(n)} \big )\qquad \text{for $i\in I$ and $n\in \N$.}
$$
Indeed, let $\mathfrak{p}\in U^\delta$; this means that $\mathfrak{p}_{\#}\in U$. Hence, let $i\in I$ such that $\mathfrak{p}_{\#}\in U_{i}$. Then, since $\mathfrak{p}_{\#}\in U_{i}$, we know that $b_{i}\notin \mathfrak{p}_{\#}$. Thus, there exists $n\in \N$ such that ${b_{i}}^{(n)}\notin \mathfrak{p}$; in other words, $\mathfrak{p}\in D\big ( {b_{i}}^{(n)} \big )$. Hence, we have found a couple $(i,n)$ such that $\mathfrak{p}\in V_{i,n}$. As families of elements, we choose
$$
\alpha_{i,n}:={a_{i}}^{(n)}\qquad \aand \qquad \beta_{i,n}:={b_{i}}^{(n)}.
$$
So, let $\mathfrak{p}\in U^\delta$ and $(i,n)$ such that $\mathfrak{p}\in V_{i,n}$. First, we have $\beta_{i,n}\notin \mathfrak{p}$. So, we have to check that
$$
t(\mathfrak{p})=\frac{\alpha_{i,n}}{\beta_{i,n}}=\frac{{a_{i}}^{(n)}}{{b_{i}}^{(n)}} \qquad \text{in $A_{\mathfrak{p}}$}
$$
By assumption, we have $\mathfrak{p}_{\#}\in U_{i}$ and so
$$
s(\mathfrak{p}_{\#}):=\frac{a_{\mathfrak{p}}}{b_{\mathfrak{p}}}=\frac{a_{i}}{b_{i}} \qquad \text{in $A_{\mathfrak{p}_{\#}}$}.
$$
Then, since these two elements have a zero derivative, and since, on the other hand, we know that $b_{\mathfrak{p}}^{(n_{\mathfrak{p}})}\notin \mathfrak{p}$ and ${b_{i}}^{(n)}\notin \mathfrak{p}$, Proposition \Ref{Constantes.Dans.Un.Anneau.Fractions} tells us that
$$
\frac{{a_{\mathfrak{p}}}^{(n_{\mathfrak{p}})}}{{b_{\mathfrak{p}}}^{(n_{\mathfrak{p}})}}=\frac{{a_{i}}^{(n)}}{{b_{i}}^{(n)}} \qquad \text{in $A_{\mathfrak{p_{\#}}}$}.
$$
Then, applying Lemma \Ref{Nulle.Dans.Traj.Nilpo.Dans.Point}.\emph{(i)} and using the fact that $A$ is reduced, we infer that
$$
\frac{{a_{\mathfrak{p}}}^{(n_{\mathfrak{p}})}}{{b_{\mathfrak{p}}}^{(n_{\mathfrak{p}})}}=\frac{{a_{i}}^{(n)}}{{b_{i}}^{(n)}} \qquad \text{in $A_{\mathfrak{p}}$}.
$$
This concludes the proof of \emph{(ii)}.

\bigskip

Now, let's indicate quickly why \emph{(i)} is true. We start with a section $s\in \mathscr{O}_{\diffSp A}^{\mathrm{(Kov)}}(U)$, a covering $(U_{i})_{i}$ of  $U$ and two families $(a_{i})_{i}$ and $(b_{i})_{i}$ of $A$ such that
$$
\fa \mathfrak{p}\in U, \quad \mathfrak{p}\in U_{i}\impl \left ( b_{i}\notin \mathfrak{p}\quad \aand \quad s(\mathfrak{p})=\frac{a_{i}}{b_{i}}\quad \text{in }A_{\mathfrak{p}} \right ).
$$
By definition, one can find Zariski open sets $W_{i}$ of $\Sp A$ such that $U_{i}=W_{i}\cap \diffSp A$, for all $i$. One can replace the $W_{i}$'s by the $\widetilde{W}_{i}$ defined by
$$
\widetilde{W}_{i}:=W_{i}\cap D(b_{i}).
$$
Then, one considers $W=\bigcup_{i} \widetilde{W}_{i}$. This is a Zariski open set of $\Sp A$ containing $U$. Let $\mathfrak{p}\in W$ and assume $\mathfrak{p}\in \widetilde{W}_{i}$ and $\mathfrak{p}\in \widetilde{W}_{j}$ for two indexes $i$ and $j$. Then,
\begin{equation}\label{egalite}
\frac{a_{i}}{b_{i}}=\frac{a_{j}}{b_{j}}\qquad \text{in }A_{\mathfrak{p}}.
\end{equation}
Indeed, $\mathfrak{p}_{\#}$ lies in $U_{i}$ and in $U_{j}$ and so
${a_{i}}/{b_{i}}={a_{j}}/{b_{j}}$ in $ A_{\mathfrak{p_{\#}}}$.
But, by Lemma \Ref{Nulle.Dans.Traj.Nilpo.Dans.Point} and for $A$ is reduced, this implies (\ref{egalite}). Then, one can define $t\in  \mathscr{O}_{\Sp A}(W)$ by setting $t(x)=a_{i}/b_{i}$ when $x\in \widetilde{W_{i}}$. The statement on unicity comes from Lemma \Ref{Nulle.Dans.Traj.Nilpo.Dans.Point}.
}

\bigskip

\section{Comparison of the Carr\`a Ferro, Keigher and Kovacic sheaves}

\subsection{Comparison of the Carr\`a Ferro and Keigher sheaves }We now come to the main result of this paper:

\smallskip

\begin{theorem}\label{Comparison.Const.Gal}
Let $X$ be a reduced $\Q$-scheme endowed with a vector field. Then, the Carr\`a Ferro sheaf and the Keigher sheaf have the same constants:
$$
\fa \, U \text{open in $X^{\vec{\mathscr{V}}}$}, \qquad C \left ( \mathscr{O}_{X^{\vec{\mathscr{V}}}}^{\mathrm{(Keigher)}}(U) \right ) \simeq C \left ( \mathscr{O}_{X^{\vec{\mathscr{V}}}}^{\mathrm{(CF)}}(U) \right ).
$$
\end{theorem}

\bigskip

\pproof{We keep the notations of the theorem. The Keigher sheaf is defined as the associated sheaf to a certain presheaf. Hence, thanks to Proposition \Ref{Associated.sheaf.constant}, the constant of the Keigher sheaf is the same as the associated sheaf to the constant of this certain presheaf. More precisely, one has
$$
C\left (  \mathscr{O}_{X^{\vec{\mathscr{V}}}}^{\mathrm{(Keigher)}}(U)\right )
\simeq
\left (  U \mapsto C\left ( \limind \limits_{
\begin{subarray}{c}
 V\text{ open in $X$} \\
 \text{and }U \subset V
\end{subarray}
}\, \mathscr{O}_{X}\left ( V \right ) \right ) \right )^\dag.
$$
Naturally, one would like to interchange $C(-)$ with $\limind $. This is not possible in general. For instance, the reader might search a example where $C\left ( A_{1} \otimes_{B} A_{2} \right )\neq C(A_{1}) \otimes_{C(B)} C\left ( A_{2} \right )$. 

\bigskip

But, in this situation, the colimit that we want to compute is of a very special kind : it is a filtered colimit. And, \emph{for such colimits}, one has
$$
C\big ( \limind_{i} A_{i} \big )\simeq \limind_{i} C(A_{i}). 
$$
Indeed, one knows, as explained  in chapters 9 and 11 of \Cite{categories}, and after Theorem 11.5.7 of the same book, that in the category $\dAnn$
filtered colimits commute with finite limits. But, given a differential ring $A$, its ring of constants can be characterized as a finite limit. More precisely, $C(A)$ can be characterized as the following kernel : 
$$
\xymatrix@C=6mm{
C(A) \ar[r]&A \ar@<1mm>[rr]^-{\Id + \partial(\cdot )\varepsilon}\ar@<-1mm>[rr]_-{\Id}&& A[\varepsilon]/\varepsilon ^2
}.
$$
Hence, $C(-)$ commutes with filtered colimits. So, one gets that
$$
C\left (  \mathscr{O}_{X^{\vec{\mathscr{V}}}}^{\mathrm{(Keigher)}}(U)\right )
\simeq
\left (  U \mapsto \limind \limits_{
\begin{subarray}{c}
 V\text{ open in $X$} \\
 \text{and }U \subset V
\end{subarray}
}\, C (\mathscr{O}_{X}\left ( V \right ))  \right )^\dag.
$$

\smallskip

Now, let $U$ be an open set of $X$. We will prove that
$$
\limind \limits_{
\begin{subarray}{c}
 V\text{ open in $X$} \\
 \text{and }U \subset V
\end{subarray}
}\, C (\mathscr{O}_{X}\left ( V \right )) = C\left ( \mathscr{O}_{X}(U^\delta) \right ).
$$
To begin with, if $V$ is a Zariski open set of $X$ that contains $U$, one has a map
$$
\varphi_{V} : C\left ( \mathscr{O}_{X}(V) \right )\longto C\left ( \mathscr{O}_{X}(U^\delta) \right ).
$$
This comes from Proposition \Ref{proposition:extension}: $\varphi_{V}$ is the composition of the extension map $\mathscr{O}_{X}(V)\longto \mathscr{O}_{X}(V^\delta)$  with the  restriction map to $\mathscr{O}_{X}(U^\delta)$. Let's prove that the $\varphi_{V}$'s make $\mathscr{O}_{X}(U^\delta)$ the sought colimit. Let $A$ be a differential ring together with compatible maps $\psi_{V} : C\left ( \mathscr{O}_{X}(V) \right )\longto A$. In particular one has a map $\psi_{U^\delta} : C \left ( \mathscr{O}_{X} ( U^\delta ) \right )\longto A$.
Let $V$ be a Zariski open set  containing $U$. All that we have to prove is that the diagram
$$
\xymatrix{
C \left (\mathscr{O}_{\mathscr{X}}  ( V  ) \right )\ar[rrrd]^{\psi_{V}}\ar[rd]_{\varphi_{V}}\\
& C \left (\mathscr{O}_{\mathscr{X}}  ( U^\delta  ) \right )\ar[rr]_{\psi_{U^\delta}}&& A
}
$$
commutes. But, in the diagram
$$
\boite{\xymatrix{
	\ar[dd]_{\mathrm{ext}_{V\to V^\delta}}C \left (\mathscr{O}_{\mathscr{X}}  ( V  ) \right )\ar[rrrd]^{\psi_{V}}\ar[rd]_{\varphi_{V}}\\
\ar@{}|{\textcircled{\scriptsize 1} \qquad}[r]& C \left (\mathscr{O}_{\mathscr{X}}  ( U^\delta  ) \right )\ar[rr]_{\psi_{U^\delta}}&& A \\
C \left ( \mathscr{O}_{\mathscr{X}} ( V^\delta ) \right ) \ar[ur]^{\text{restric}} \ar[rrru]_{\psi_{V^\delta}}\ar@{} ^{\qquad\textcircled{\scriptsize 2} \qquad}[urrr]
}},
$$
the diagram $\textcircled{\scriptsize 1}$ commutes by definition of $\varphi_{V}$, the diagram $\textcircled{\scriptsize 2}$ commutes for the  $\psi_{W}$'s form a compatible family, and the big triangle
$$
\boite{\xymatrix@C=1.3cm{
	\ar[dd]_{\mathrm{ext}_{V\to V^\delta}}C \left (\mathscr{O}_{\mathscr{X}}  ( V  ) \right )\ar[rrrd]^{\psi_{V}}\\
&&&A\\
C\left ( \mathscr{O}_{\mathscr{X}} ( V^\delta )  \right ) \ar[rrru]_{\psi_{V^\delta}}
}}
$$
commutes for the same reason, and for the restriction map and the extension map are inverse one of each other. So, the last triangle
$$
\boite{\xymatrix{
C \left (\mathscr{O}_{\mathscr{X}}  ( V  ) \right )\ar[rrrd]^{\psi_{V}}\ar[rd]_{\varphi_{V}}\\
& C \left (\mathscr{O}_{\mathscr{X}}  ( U^\delta  ) \right )\ar[rr]_{\psi_{U^\delta}}&& A
}}
$$
indeed commutes and $C\left ( \mathscr{O}_{X}(U^\delta) \right )$ is the colimit that we wanted to compute.

\bigskip

Now, the end of the proof is easy. Since the constant of a sheaf is still a sheaf, the presheaf
$$
U \mapsto C\left ( \mathscr{O}_{X}(U^\delta) \right ),
$$
which is actually the constant of the Carr\`a Ferro sheaf, is a sheaf. So, it is its own associated sheaf.
}

\medskip

\remark{In general, the sheaves $\mathscr{O}_{X^{\vec{\mathscr{V}}}}^{\mathrm{(Keigher)}}$ and $ \mathscr{O}_{X^{\vec{\mathscr{V}}}}^{\mathrm{(CF)}}$ are not isomorphic. For instance, if $\mathscr{X}$ is $\A{1}{\C}$ with the constant vector field, as we already told, $X^{\vec{\mathscr{V}}}$ contains only one element, the generic point. The global sections are $\C[t]$ for the Carr\`a Ferro sheaf and $\C(t)$ for the Keigher sheaf.
}

\bigskip

\subsection{Comparison of the Keigher sheaf and the Kovacic sheaf}
We now prove the following proposition, that compares the two classical sheaves over $\diffSp A$:

\begin{proposition}\label{Compare.Sheaf.Aff}
Let $A$ be a differential ring. Then,
$$
\mathscr{O}_{\diffSp A}^{\mathrm{(Keigher)}}\simeq \mathscr{O}_{\diffSp A}^{\mathrm{(Kov)}}.
$$
\end{proposition}

\smallskip

As an immediate consequence of Theorem \Ref{Comparison.Const.Gal} and of the latter proposition, one gets:

\smallskip

\begin{corollary}
Let $A$ be $\Q$-reduced differential algebra. Then, 
$$
C \big ( \mathscr{O}_{\diffSp A}^{\mathrm{(Keigher)}} \big )\simeq C \big ( \mathscr{O}_{\diffSp A}^{\mathrm{(Kov)}} \big ) \simeq C \big ( \mathscr{O}_{\diffSp A}^{\mathrm{(CF)}} \big ).
$$
\end{corollary}

\smallskip

\pproofbis{of Proposition \Ref{Compare.Sheaf.Aff}}
{First, let us remark\footnote{It follows, for the Keigher sheaf, from the fact that a sheaf and its restriction to a subset have the same stalks.}  that $\mathscr{O}^{\textrm{(Kov)}}_{\diffSp A}$ and $\mathscr{O}^{\textrm{(Keigher)}}_{\diffSp A}$ have the same stalks: for all $\mathfrak{p}\in \diffSp A$, the stalks at $\mathfrak{p}$ are isomorphic to $A_{\mathfrak{p}}$. So, to prove that they are isomorphic, it is sufficient to show that there exists a morphism between them, inducing isomorphisms on stalks. Let us indicate how to construct such a morphism $\mathscr{O}^{\textrm{(Keigher)}}_{\diffSp A} \longto \mathscr{O}^{\textrm{(Kov)}}_{\diffSp A}$. By the universal property of the associated sheaf, it is sufficient to built a similar morphism
$$
\mathop{\limind}\limits_{\begin{subarray}{c}
U\supset V \\[1mm]
V\text{ Zariski open} \\ \text{in }\Sp A
\end{subarray}}  \mathscr{O}_{\Sp A}\left ( V \right ) \longto \mathscr{O}^{\textrm{(Kov)}}_{\diffSp A} (U),
$$
functorial in $U$. But, to define such a morphism, it is sufficient to consider a compatible family of morphisms
$$
\mathscr{O}_{\Sp A}(V) \longto \mathscr{O}^{\textrm{(Kov)}}_{\diffSp A} (U)
$$
for all Zariski open set $V$ containing $U$. If one consider the Hartshorne-like definition of $\mathscr{O}_{\Sp A}$, it is easy to define these maps, by restriction.}

\bigskip

\appendix

\section{The associated sheaf in the differential \hyphenation{context}context}

The two goals of this appendix are:
\begin{itemize}
\item[\emph{(i)}] to explain why the existence of the associated sheaf $\mathscr{F}^{\dagger}$ in the context of sheaves and presheaves of sets implies its existence in the context of differential rings;
\item[\emph{(ii)}] to explain why the functor $\mathscr{F}\mapsto \mathscr{F}^{\dagger}$ commutes with the functor of constants.
\end{itemize} 

\subsection{Statement of the first result} To begin with, we fix some notations. If $X$ is a topological space, we denote by 
\begin{align*}
\PrSh{X}&&\AbPrSh{X} && \RngPrSh{X}& &  \dRngPrSh{X} \\
\Sh{X}&&\AbSh{X} && \RngSh{X}& &  \dRngSh{X}
\end{align*}
the respective categories of presheaves and sheaves of sets, abelian groups, rings and differential rings. These categories come naturally with forgetful functors. We also denote by 
\begin{align*}
&C_{(\textit{Sh})} : \dRngPrSh{X}\longto \RngPrSh{X} \\ \aand \qquad & C_{(\textit{PrSh})} : \dRngSh{X} \longto \RngSh{X}
\end{align*}
the functors that associate to a (pre)sheaf $\mathscr{F}$ of differential rings the (pre)sheaf of rings $U \mapsto C\left ( \mathscr{F}(U) \right )$.
Finally, we also recall that one denotes
$$
\xymatrix@C=1.5cm{
\mathscr{C}\ar@{}|-{\top}[r]\ar@<1.7mm>[r]^{F} & \mathscr{D} 
\ar@<1.7mm>[l]^{G}
}
$$
when $(F,G)$ establishes an adjunction between $\mathscr{C}$ and $\mathscr{D}$, \ie when $G$ is  left adjoint to $F$. We want to prove:\smallskip

\begin{theorem}
Let $X$ be a topological space and let $\mathscr{F}$ be a presheaf of sets over $X$. Then, the sheaf of sets $\mathscr{F}^\dag$, associated to $\mathscr{F}$, can be endowed with a structure of sheaf of abelian groups (resp. rings, differential rings) when $\mathscr{F}$ has a structure of presheaf of abelian groups (resp. rings, differential rings) .
\end{theorem}

 But, more precisely, what we will prove is the following

\begin{theorem}\label{theo.asso.sheaf.diff.case}
Given a topological space $X$, there exist four adjonctions
\begin{align*}&
  \xymatrix@C=1.1cm{
\Sh{X}\ar@{}|-{\top}[r]\ar@<1.7mm>[r]^-{\omega} & \PrSh{X} \ar@<1.7mm>[l]^-{*}
} &&
\xymatrix@C=1.1cm{
\AbSh{X}\ar@{}|-{\top}[r]\ar@<1.7mm>[r]^-{\omega_{\textit{Ab}}} & \AbPrSh{X} \ar@<1.7mm>[l]^-{*_{\textit{Ab}}}}
\\&
\xymatrix@C=1.1cm{
\RngSh{X}\ar@{}|-{\top}[r]\ar@<1.7mm>[r]^-{\omega_{\textit{Rng}}} &\RngPrSh{X} \ar@<1.7mm>[l]^-{*_{\textit{Rng}}}}
& &
\xymatrix@C=1.1cm{
\dRngSh{X}\ar@{}|-{\top}[r]\ar@<1.7mm>[r]^-{\omega_{\textit{Rng}^\partial}} & \dRngPrSh{X} \ar@<1.7mm>[l]^-{*_{\textit{Rng}^\partial}}}
\end{align*}
making the following diagram commute:
$$
\xymatrix@C=1.2cm@R=1.2cm{
\ar@<1.5mm>[d]^-{\omega_{\textit{Rng}^\partial}}	\dRngSh{X} \ar[r]_-{\omega_{1}} &\ar@<1.5mm>[d]^-{\omega_{\textit{Rng}}}	\RngSh{X} \ar[r]_-{\omega_{2}} &\ar@<1.5mm>[d]^-{\omega_{\textit{Ab}}}	 \AbSh{X}\ar[r]_-{\omega_{3}} &	 \Sh{X} \ar@<1.5mm>[d]^-{\omega}	 \\
\ar@<1.5mm>[u]^-{*_{\textit{\textit{Rng}}^\partial}}\dRngPrSh{X} \ar[r]_-{\omega_{4}} &\ar@<1.5mm>[u]^-{*_{\textit{Rng}}}\RngPrSh{X} \ar[r]_-{\omega_{5}} &\ar@<1.5mm>[u]^-{*_{\textit{Ab}}} \AbPrSh{X} \ar[r]_-{\omega_{6}} &\ar@<1.5mm>[u]^-{*} \PrSh{X} \\
}.
$$
Moreover, the functors $*$, $*_{\textit{Ab}}$, $*_{\textit{Rng}}$ and $*_{\textit{Rng}^\partial}$ commute with general colimits and finite limits.
\end{theorem}

\smallskip

With these notations, the left adjoint functor to $\omega:\Sh{X}\to \PrSh{X}$, denoted here by $* : \PrSh{X}\to \Sh{X}$, is the functor that associates to a presheaf $\mathscr{F}$ its associated sheaf $\mathscr{F}^\dag$.

\smallskip

\subsection{Proof of Theorem \Ref{theo.asso.sheaf.diff.case}}
To begin with, we assume that the existence of the associated sheaf, and the fact that it commutes with finite limits, is known for presheaves of sets. It is proved, for instance, in \Cite{hartshorne} or \Cite{sheavestopos}. We denote by $* : \PrSh{X}\to \Sh{X}$ the functor that maps $\mathscr{F}$ to its associated sheaf $\mathscr{F}^\dag$.

\bigskip

Now, let $\mathscr{C}$ be a category with finite products. We denote by $\mathbf{1}$ a terminal object of $\mathscr{C}$. Following \Cite[Ch. II, \S 7]{sheavestopos}, but the interested reader should also consult \Cite[Ch. 11]{categories}, we consider the category $\mathbf{Ab}\left ( \mathscr{C} \right )$ of \emph{abelian group objects of $\mathscr{C}$}. It is defined as follows :
\begin{itemize}
\item[---] the objets of $\mathbf{Ab}\left ( \mathscr{C} \right )$ are $4$-uplets $(X, a, v, u)$ where $X\in \ob{\mathscr{C}}$ and where $a : X \times X \to X$, $v : X\to X$ and $u : \mathbf{1}\to X$ are arrows satisfying some conditions. Intuitively, one wants $a$ to be the addition law, $v$ to be the subtraction law and $u$ to be the zero.
\item[---] the arrows of $\mathbf{Ab} \left ( \mathscr{C} \right )$ are arrows $f : X \longto Y$ that commutes with addition.
\end{itemize}
For instance, $\mathbf{Ab}\left ( \mathbf{Sets} \right )$ is isomorphic, as a category, to the category of abelian groups. Similarly, for every topological space $X$, the categories $\mathbf{Ab}\left ( \PrSh{X}\right )$ and $\mathbf{Ab}\left ( \Sh{X}\right )$ are isomorphic to $\AbPrSh{X}$ and to $\AbSh{X}$.
\bigskip

Now, let $X$ be a topological space and let $\mathscr{F}$ be presheaf of abelian groups. We want to construst the associated sheaf $\mathscr{F}^\dag$ to $\mathscr{F}$. First, one can see $\mathscr{F}$ as an object of $\mathbf{Ab}\left ( \PrSh{X}\right )$: $\mathscr{F}$ is preasheaf of \emph{sets} endowed with maps 
$$
a : \mathscr{F}\times \mathscr{F}\longto \mathscr{F},\qquad v : \mathscr{F} \longto \mathscr{F}\quad \aand\quad u : \left \{\star\right \} \longto \mathscr{F}
$$
where $\left \{\star\right \} $ denotes the final object of $\PrSh{X}$. Then, one can apply to $\mathscr{F}$ and to these maps the functor $*$. Since, $*$ commutes with finite limits, one gets
$$
*(a) : \mathscr{F}^\dag \times \mathscr{F}^\dag \longto  \mathscr{F}^\dag,
\qquad *(v) : \mathscr{F}^\dag \longto \mathscr{F}^\dag\quad \aand\quad *(u) : \left \{\star\right \}^\dag \longto \mathscr{F}^\dag.
$$
Furthermore, these maps still verify the axioms of $\mathbf{Ab}\left ( \mathscr{C} \right )$, since $*$ maps commutative diagrams to commutative diagrams. Therefore, $\mathscr{F}^\dag$ has naturally a structure of sheaf of abelian groups. One can also verify that $*$ maps additive morphisms of additive morphisms. Thus, one has a functor
$$
*_{\textit{Ab}} : \AbPrSh{X}\longto \AbSh{X}
$$ 
and one can prove that it is the left adjoint that we were looking for. Last, since $*_{\textit{Ab}}$ is left adjoint to $\omega_{\textit{Ab}}$, one knows that it commutes with all colimits. For finite limits, one proceeds as follows:
\begin{enumerate}
\item Let $(\mathscr{F}_{i}, \varphi_{ij})$ be a finite system of abelian presheaves and let $\varphi_{i} : \mathscr{F}\longto \mathscr{F}_{i}$ be its limit in $\AbPrSh{X}$. Then, $\mathscr{F}$, seen as a presheaf of \emph{sets} is still a limit. This comes, for instance, from the fact that the functor $\omega_{6} : \AbPrSh{X}\longto \PrSh{X}$ has a left adjoint. This adjoint maps a presheaf of sets $\mathscr{G}$ to the presheaf of abelian groups $U \mapsto \Z^{(\mathscr{G}(U))}$.
\item Hence, $\varphi_{i} :\mathscr{F} \longto \mathscr{F}_{i}$ is still a limit, seen in $\PrSh{X}$. For $*$ commutes with finite limits, one gets that $*(\varphi_{i}) : {\mathscr{F}}^\dag\longto {\mathscr{F}_{i}}^\dag$ is a limit in $\Sh{X}$. Furthermore, by definition, ${\mathscr{F}_{i}}^\dag$, $\mathscr{F}^\dag$ can be seen as sheaves of abelian groups and the morphism $*(\varphi_{i})$ are additive.
\item Last, if $\mathscr{G}$ is a sheaf of abelian groups given with morphisms $\psi_{i} : \mathscr{G} \longto {\mathscr{F}_{i}}^\dag$, one wants to find an arrow $f : \mathscr{G} \longto {\mathscr{F}}^\dag$ that factorizes the $\psi_{i}$. For ${\mathscr{F}}^\dag$ is a limit in $\Sh{X}$, one can find such an arrow $f$, in $\Sh{X}$. But then, one has to prove that this arrow is additive. This comes from the additiveness of the $\psi_{i}$ and the unicity of factorizations.
\end{enumerate}
\bigskip

So, this is how one can prove the existence of the left adjoint $*_{\textit{Ab}}$ and its properties. The same arguments apply for (pre)sheaves of rings and differential rings: one remarks that it is possible to defines ring objects and differential ring objets in a category $\mathscr{C}$ with finite limits, and that this definition only involves finite products, maps and commutative diagrams.

\subsection{Associated sheaves and constants }Now, we prove
\begin{proposition}\label{Associated.sheaf.constant}
Let $X$ be a topological space. Then, the diagram
$$
\xymatrix@C=1.5cm@R=1cm{
\dRngSh{X} \ar[r]_-{C_{(\textit{Sh})} } &\RngSh{X}  \\
\ar[u]^-{*_{\textit{Rng}^\partial}}\dRngPrSh{X} \ar[r]_-{C_{(\textit{PrSh})} } &\ar[u]^-{*_{\textit{Rng}}}\RngPrSh{X} \\
}
$$
commutes, up to isomorphism.
\end{proposition}

\smallskip

This means that, if $\mathscr{F}$ is a presheaf of differential rings, when one wants to compute the constant of $\mathscr{F}^\dag$, it suffices to compute the associated sheaf to $C\left ( \mathscr{F} \right )$: in other words, $C ( \mathscr{F} )^\dag\simeq C ( \mathscr{F}^\dag )$.

\smallskip

\pproof{Let $X$ be a topological space and let $\mathscr{F}\in \dRngPrSh{X}$. The constant of $\mathscr{F}$, denoted by $C\left ( \mathscr{F} \right )$ is a  finite limit; more precisely, it is a kernel (in $\PrSh{X}$, and in $\Sh{X}$ for sheaves):
$$
\xymatrix{
C\left ( \mathscr{F} \right )\ar[r] &\mathscr{F}\ar@<1.5mm>[r]^-\partial \ar@<-1.5mm>[r]_-0 &\mathscr{F}.
}
$$
For $*$ commutes with finite limits, one has
$$
\xymatrix{
C\left ( \mathscr{F} \right )^\dag \ar[r] &\mathscr{F}^\dag \ar@<1.5mm>[r]^-\partial \ar@<-1.5mm>[r]_-0 &\mathscr{F}^\dag:
}
$$
hence, $C\left ( \mathscr{F} \right )^\dag$ is isomorphic to $C\left ( \mathscr{F}^\dag \right )$, as sheaves of sets. But, this enough to infer that they are isomorphic as sheaves of rings. 
Indeed, 
first, the map $C\left ( \mathscr{F} \right )^\dag \longto \mathscr{F}^\dag$ is a morphism of sheaves of differential rings; second, if $U$ is any open set, the map $C\left ( \mathscr{F} \right )^\dag(U) \longto \mathscr{F}^\dag(U)$ is injective (this is because it is isomorphic to $C ( \mathscr{F}^\dag )(U) \longto \mathscr{F}^\dag$); third, as sets, 
$C\left ( \mathscr{F}^\dag \right )(U) $ and $C\left ( \mathscr{F} \right )^\dag(U)$ are isomorphic.
}

\newpage
\bibliography{Schemes_and_vector_fields}

\end{document}